\documentclass[11pt]{amsart}

\usepackage{graphicx}% Required for inserting images
\usepackage{amsmath,amssymb,amsthm}
\usepackage[colorlinks,backref,citecolor=myGreen,linkcolor=darkred]{hyperref}
\usepackage{thmtools, thm-restate}
\usepackage{cleveref}
\usepackage{imakeidx}
\usepackage{algpseudocode}
\usepackage{algorithm}
\usepackage{mathtools}
\usepackage{subcaption}
\usepackage{caption}
\usepackage{tikz}
\usetikzlibrary{cd}
\usepackage{xcolor}
\usepackage{soul}
\usepackage{faktor}
\usepackage{mathscinet}
\usepackage{cancel}
\usepackage{longtable}
\usepackage{booktabs}
\usepackage{enumerate}
\usepackage[mode=buildnew]{standalone}
\usepackage{xr}
\usepackage{bm}

\usetikzlibrary{positioning}
\allowdisplaybreaks

\usepackage[colorinlistoftodos, bordercolor=black, backgroundcolor=orange!20, linecolor=orange, textsize=scriptsize]{todonotes}

\makeindex[intoc,name=general]

\newcommand\DEF[1]{\emph{#1}}

\newcommand{\gitrepo}{\text{\url{https://github.com/dmg-lab/GroebnerQuantumSym.jl}}}

\definecolor{myGreen}{rgb}{0.18039216 0.49803922 0.09411765}
\definecolor{darkred}{rgb}{0.65, 0.0, 0.0}

\newcommand{\C}{\mathbb{C}}

\newcommand{\N}{\mathbb{N}}

\newcommand{\Z}{\mathbb{Z}}

\newcommand{\hlmath}[1]{#1}

\newcommand\rs{\mathsf{rs}}
\newcommand\cs{\mathsf{cs}}
\newcommand\ip{\mathsf{ip}}
\newcommand\inj{\mathsf{inj}}
\newcommand\wel{\mathsf{wel}}
\newcommand\rwel{\mathsf{rwel}}
\newcommand\rinj{\mathsf{rinj}}

\newcommand\rrs{\mathsf{rrs}}
\newcommand\rcs{\mathsf{rcs}}
\newcommand\rinjcs{\mathsf{rinjcs}}
\newcommand\rwelcs{\mathsf{rwelcs}}

\DeclareMathOperator{\Span}{span}
\DeclareMathOperator{\lm}{lm}
\DeclareMathOperator{\len}{len}
\DeclareMathOperator{\lc}{lc}
\DeclareMathOperator{\lt}{lt}
\DeclareMathOperator{\im}{im}

\newcommand{\fL}{\mathfrak{L}}

\newcommand{\QSn}[1]{\mathfrak{S}_{#1}}

\newcommand{\bg}[1]{\mathsf{bg}^{(#1)}}
\newcommand{\rbg}[1]{\mathsf{rbg}^{(#1)}}

\newcommand{\MathInHeader}[1]{\texorpdfstring{$#1$}{TEXT}}

\newtheorem{theorem}{Theorem}[section]
\theoremstyle{plain}

\newtheorem{corollary}[theorem]{Corollary}

\newtheorem{example}[theorem]{Example}

\newtheorem{lemma}[theorem]{Lemma}
\newtheorem{definition}[theorem]{Definition}

\newtheorem{proposition}[theorem]{Proposition}
\newtheorem{remark}[theorem]{Remark}

\title{gb4quantumGroups}
\author{}
\date{May 2024}

\begin{document} 
%\tableofcontents
\bibliographystyle{plain}

\title{Finite Gröbner bases for quantum symmetric groups}

\author[L. Schmitz]{Leonard Schmitz}
\address{
Technische Universität Berlin, Algebraic and Geometric Methods in Data Analysis}
\email{lschmitz@math.tu-berlin.de}

\author[M. Wack]{Marcel Wack}
\address{
Technische Universität Berlin, Chair of Discrete Mathematics / Geometry}
\email{wack@math.tu-berlin.de}


\date{March 2025}

%\listoftodos{}

\begin{abstract}Non-commutative Gröbner bases of two-sided ideals are not necessarily finite. Motivated by this, we provide a closed-form description of a finite and reduced Gröbner bases for the two-sided ideal used in the construction of Wang's quantum symmetric group. In particular, this proves that the word problem for quantum symmetric groups is decidable. 
\end{abstract}

\maketitle

\section{Introduction}

Gr\"obner bases were introduced  by Bruno Buchberger in his seminal  Ph.D. thesis \cite{buchberger1965algorithmus}, and have since revolutionized various  disciplines,  including computational algebraic geometry \cite{cox1997ideals,sturmfels2002solving}, commutative algebra \cite{eisenbud2013commutative},  discrete geometry \cite{joswig2013polyhedral,sturmfels1996grobner}, algebraic statistics \cite{sullivant2023algebraic}, and many others, e.g. \cite{michalek2021invitation}. A Gr\"obner basis is a particular set of generators for a given ideal in a commutative polynomial ring. Every
set of polynomials can be transformed into a Gr\"obner basis via Buchberger's algorithm.
\par
Several extensions to more general algebraic structures have been developed, for example Gr\"obner bases of modules \cite{moller1986new,kreuzer2000computational}, extensions over  principal ideal domains \cite{kandri1988computing}, Ore algebras \cite{KANDRIRODY19901} and many others.  
\par
This work is rooted in the non-commutative extension of Gr\"obner bases  for free associative algebras \cite{mora1988groebner,bergman1978diamond}. In contrast to Buchberger's setting over commutative polynomials, there does not exist a finite Gr\"obner basis for an arbitrary two-sided ideal.  
Closely related to this, the \DEF{word problem in free associative algebras} is the problem of deciding whether two given polynomials are equivalent
modulo a given two-sided ideal.
The equivalent  ideal membership problem  is only \emph{semi-decidable} \cite{pritchard1996ideal},  meaning that there is a suitable procedure which terminates if and only if the ideal membership takes place, or runs forever otherwise. For a semi-decision procedure based on the so-called \emph{letterplace embedding},  compare \cite{LASCALA20091374}. Nonetheless, if there exists a finite Gr\"obner basis for a given ideal, then the associated word problem becomes decidable. 
\par 
Several applications have recently emerged for non-commutative Gr\"obner bases, e.g. formalization techniques for matrix identities \cite{HRR19,hofstadler2024universaltruthoperatorstatements,SL20} or symmetry decision processes for quantum automorphism groups of graphs \cite{Lev22} and matroids \cite{corey2023quantum}. The latter is the motivation for this work.
\par
   We prove that the word problem is decidable for every quantum symmetric group $\mathfrak{S}_n$ of \DEF{size} $n\in\N$. 
  For this consider the notion of quantum symmetric groups introduced by Wang, which in turn fits into Woronowicz's theory of compact quantum groups \cite{woronowicz1987compact}.
Wang in his work in 1998 \cite{wang1998quantum} found that these quantum groups give the possibility to quantum permute the elements of a finite set, extending the classical symmetric group and thus creating the quantum symmetric group roughly as the quotient of a $C^*$-algebra by some ideal $I_n$ (equipped with some additional structure).
Since the inception of the quantum symmetric group, the field has developed to define as subgroups of the quantum symmetric group the quantum automorphism group of graphs \cite{bichon2003quantum} and other combinatorial objects. 
A large part of the research in these fields is concerned with checking whether a given structure possesses quantum symmetries in addition to its classical symmetries. 
Since subgroups of the quantum symmetric group are quotients over some ideal,  
searching for quantum symmetries amounts to solving an ideal membership problem for the defining ideal, which motivates the search for a finite Gr\"obner basis of the latter. 
This has been done computationally for some examples in \cite{Lev22}. 
In general it is conjectured to be undecidable whether a given graph is quantum symmetric or not, making it all the more remarkable that we can show that the word problem is decidable for any quantum symmetric group.
For this we translate the problem into non-commutative factor algebras and provide a finite closed-form description of the Gr\"obner basis $G_n$ for the underlying two-sided ideal of quantum symmetries $I_n$, and hence a linear decision procedure for word problems in quantum symmetric groups. 
\par
Even the existence of such a closed-form Gr\"obner bases for an entire family of non-trivial ideals is remarkable and rare. One of the most prominent scenarios are (commutative)  Pl\"ucker relations  \cite[Theorem 5.8]{michalek2021invitation} and their associated ideals. In non-commutative algebras, however, the existence of these closed-form Gr\"obner bases is even more essential since this solves the otherwise undecidable word problem for the entire class of resulting factor algebras. 
\par

In constructing the Gröbner basis we are faced with rather large parameterised calculations for which even writing down the equation would be a challenge in itself. 
As a solution, we offer a rather novel approach to perform large parameterised (non-commutative) polynomial calculations using the computer algebra system \texttt{OSCAR} \cite{OSCAR}. 
A detailed explanation is given in \Cref{sec:predicates}, and the necessary computations are available at
\[
  \gitrepo
\]

This manuscript is organized into the following parts: in \Cref{sec:mainresult} we provide all elementary definitions in order to state our main result, \Cref{thm:finite_gb}, without going into technical details yet. \Cref{sec:quantumSymGroup}
is an introduction to the quantum symmetric group. In particular, it covers the important notions \emph{transposition} and \emph{symmetry}. 
In \Cref{sec:primGBs} we give a primer on Gr\"obner bases to finally prove our main result, \Cref{sec:main_result}.
This central proof section is organized into three subsections: \Cref{subsec:proof_first_interreduce} constructs an interreduced generating set for our two-sided ideal $I_n$, \Cref{subsec:mainproof_overlaps_bg} covers the corresponding Gr\"obner basis $G_n$, and in \Cref{subsec:proof_veryBadGuys} we cover the remaining overlap relations that have not been considered by the previous parts.
Finally, in \Cref{sec:generalProof} we conclude with our general computational proof strategy for the larger and in \Cref{sec:main_result} omitted Gr\"obner certificates. 
\subsection{The main result}\label{sec:mainresult}

Following the notation in  \cite{corey2023quantum,Lev22} let
$$\mathbb{C}\langle n^2\rangle:=\mathbb{C}\langle u_{ij}\mid 1\leq i,j\leq n\rangle$$ 
be the  \DEF{free associative algebra} over $n^2$ \DEF{symbols} $u_{ij}$. 
For any $1\leq i,j,k\leq n$ with $k\not=j$ we denote by 
\begin{align}
\label{eq:def_rscs}
&\rs_i:=\sum_{1\leq \alpha\leq n}u_{i\alpha}-1,&&\cs_i:=\sum_{1\leq \alpha\leq n}u_{\alpha i}-1,\\
\label{eq:def_inj_rwel}
&\inj_{jik}:=u_{ji}u_{ki},&&\wel_{ijk}:=u_{ik}u_{ij},\\
\label{eq:def_ip}
&\ip_{ij}:=u_{ij}^2-u_{ij}
\end{align}  the \DEF{row} and  \DEF{column sum relations} \eqref{eq:def_rscs}, 
 the \DEF{orthogonal relations} \eqref{eq:def_inj_rwel}, and the \DEF{idempotent relations} \eqref{eq:def_ip},  
 respectively. These relations 
define the algebraic 
\DEF{quantum symmetric group} 
\begin{equation*}
  \QSn{n} := \faktor{\C\langle n^2 \rangle}{I_n}
\end{equation*}
as a factor algebra modulo the two-sided \DEF{ideal of quantum symmetries}, 
$$I_n:=\left\langle \rs_i,\cs_i,\ip_{ij}, \inj_{jik},\wel_{ijk}\;\begin{array}{|l}\;1\leq i,j,k\leq n\\\text{ with }j\not=k\end{array}\right\rangle.
$$
Note that $I_n$ is finitely generated by $2n(n^2+1)$ inhomogeneous generators. 
We introduce two other families for which it will become clear later that these are certain orthogonal relations reduced modulo the generators of $I_n$. 
\begin{definition}\label{def:rinjRwel}
 For all $n\geq 4$, the ideal $I_n$ contains the \DEF{reduced orthogonal relations}
\begin{align*}
  \rinj_{jk}&:= \sum_{3\leq\alpha \leq n} u_{j2}u_{k\alpha} - \sum_{3\leq\alpha \leq n} u_{j\alpha}u_{k1} + u_{k1} - u_{j2} \\
  \rwel_{jk}&:= \sum_{3\leq\alpha \leq n} u_{2j}u_{\alpha k} - \sum_{3\leq\alpha \leq n} u_{\alpha j}u_{1k} + u_{1k} - u_{2j}
\end{align*}
where $2\leq j,k\leq n$ with $j\not=k$.
\end{definition}

Using these additional relations, we can state our main result. 
For readability, we omit the conditions for the indices when the resulting relations are not defined. 
%and we are inteste%\marcel{what is a good name for those?}

\begin{theorem}\label{thm:finite_gb}
For all $n\geq 4$ the ideal $I_n$ has a finite, monic and reduced Gröbner basis 
\begin{align*}
G_n:=\{\cs_1\}
&\cup\left\{\begin{array}{l|}\cs_i,\rs_i,\ip_{ij},\inj_{ijk}\\\wel_{ijk},\mathsf{rinj}_{kj},\mathsf{rwel}_{kj}\end{array}\; \; i,j,k\not=1\right\}
\\
&\cup\left\{u_{k2} \inj_{j3i}-\rinj_{kj}u_{i3}\;\;\begin{array}{|l}i,j,k\not=1\text{ and}\\(k,j)\not=(2,3)\not=(j,i)\end{array}\right\}
\\
&\cup\left\{u_{2k}\wel_{3ji}-\rwel_{kj}u_{3i}\;\;\begin{array}{|l} i,j,k\not=1,\\(k,j)\not=(2,3)\not=(j,i)\\
\text{and }(k,j,i)\not=(2,4,3)\end{array}\right\}
\end{align*}
with respect to the graded lexicographic order via row-wise ordering in $(u_{ij})_{1\leq i,j\leq n}$. Its explicit cardinality is given by the cubic polynomial $\# G_n=4n^3 -15n^2 +16n-2$. 
\end{theorem}

% Its  leading monomials are 
%$$\lm(G^{0})\cup\left\{u_{k2}u_{j4}u_{i3},u_{2k}u_{4j}u_{3i}\; \begin{array}{|l}k,j,i\in[2,n]\text{ where}\\
%(k,j)\not=(2,3)\not=(j,i)\end{array}\right\}$$
%\leo{i would sugest to omit the leadning monomilas, but to give the exact Gr\"oebner basis and the number of elements instead.}

\begin{corollary}
The word problem in $\QSn{n}$ is decidable. 
\end{corollary}

Before continuing with the proof of \Cref{thm:finite_gb}, first connect it to the notion of a quantum symmetric group.

\section{Quantum symmetric groups}\label{sec:quantumSymGroup}
While the theory of quantum groups is extensive and deeply rooted in the theory of functional analysis, we will only give a brief overview of the reasons for studying the algebraic core of the quantum symmetric group, so that no further knowledge is required. The specific type of quantum group we need to know about are compact matrix quantum groups, introduced by Woronowicz \cite{woronowicz1987compact,woronowicz1991remark} in 1987.
\begin{definition}\label{def:Woronowicz}
  A \emph{compact matrix quantum group} $G$ is a pair $(C(G),u)$, where $C(G)$ is a unital $C^*$-algebra that is  generated by the $C(G)$-valued entries of the $n\times n$ matrix $$u := (u_{ij})_{1\leq i,j \leq n}.$$
  Furthermore, $u$ and its matrix transpose $u^{\top}$ must be invertible in $C(G)^{n\times n}$, and 
\begin{align*}
  \Delta : C(G) &\rightarrow C(G) \otimes C(G), \\
  u_{ij} &\mapsto \sum_{k=1}^n u_{ik} \otimes u_{kj}
\end{align*}
has to be a $^*$-homomorphism. The map $\Delta$ is called the \emph{coproduct} of $G$.
\end{definition}
To be even more specific, we will only talk about quantum symmetric groups developed by Wang \cite{wang1998quantum} in 1998.
\begin{definition}\label{def:wang}
  The \emph{quantum symmetric group} $S_n^+ := (C(S_n^+),u)$ is the compact matrix quantum group given by a matrix $u := (u_{ij})_{1\leq i,j \leq n} \in C(S_n^+)$ subject to
  \begin{align*}
    u_{ij}^* &= u_{ij} = u_{ij}^2, &1\leq i,j \leq n\\
    \sum_{i=1}^n u_{ij} &= \sum_{k=1}^n u_{kh}  =1  &1\leq j,h \leq n.
  \end{align*}
Together with the universal $C^*$-algebra $C(S_n^+)$ generated by entries of $u$.
\end{definition}
Let us restate the following well known fact, e.g. \cite{Manvcinska2020nonlocalgames}.

\begin{theorem}
  The $C^*$-algebra $C(S_n^+)$ is commutative if $n< 4$, i.e. $u_{ij}u_{kl} = u_{kl}u_{ij}$ for all $1\leq i,j,k,l \leq n$.
  If $n\geq 4$, then $C(S_n^+)$ is non-commutative.
\end{theorem}
This already explains the assumption $n\geq 4$ for our main theorem: existence of finite Gr\"obner bases is well-known for commutative algebras, e.g. \cite{cox1997ideals}. In other words, we could trivially extend \Cref{thm:finite_gb} in the remaining cases $1\leq n <4$. 

The notion of a compact matrix quantum group according to \Cref{def:Woronowicz,def:wang} requires the language of $C^*$-algebras, while in \Cref{sec:mainresult} we used the language of free associative algebras for our introduction of the quantum symmetric group. This is justified since compact matrix quantum groups admit a dense $^*$-subalgebra generated by the entries of $u$. For a detailed explanation, we refer the reader to \cite[Section 5]{timmermann2008invitation}.

The second difference is the absence of the \emph{orthogonal} relations in the definitions from above.
This is motivated by  \cite{raeburn2005graph}
where it is shown that  all orthogonal relations can already be implied by those based on \emph{projection/idempotence} and \emph{row/column sums}. 
Compare also \cite{speicher2019quantum} for a nore detailed explanation of this is in the more general setting of graph automorphism groups.
Having said this, we can proceed with the definition of the algebraic version used in the introduction, \Cref{sec:mainresult}. 
\begin{definition}
  Let $n \in \N$. The \emph{algebraic quantum symmetric group}  
  \begin{equation}\label{eq:def_quantumGroup}
    \QSn{n} := \faktor{\C\langle n^2 \rangle}{I_n}    
  \end{equation}
  is the free associative algebra over $n^2$ symbols $u_{ij}$ modulo the two-sided ideal of \emph{quantum symmetries} $I_n$
 % $$I_n := \langle \rs_i,\cs_i,\inj_{ijk},\wel_{ijk},\ip_{ij} \rangle$$
  generated by the relations \eqref{eq:def_rscs}, \eqref{eq:def_inj_rwel} and \eqref{eq:def_ip}. 
  %\begin{align*}
  %  \rs_i &:= \sum_{j=1}^nu_{ij}-1, &\cs_i := \sum_{i=1}^nu_{ij}-1 \quad &1\leq i,j \leq n\\
   % \inj_{ijk} &:= u_{ij}u_{kj} \quad (k \neq i) &\wel_{ijk} := u_{ij}u_{ik} \; (k \neq j) \quad  &1\leq i,j,k \leq n\\
   % \ip_{ij} &:= u_{ij}^2 - u_{ij}, &&1\leq i,j \leq n.
  %\end{align*}
\end{definition}
To lift this definition to be the dense $^*$-subalgebra of the quantum symmetric group, one could equip it with the $^*$-involution $(u_{ij})^* = u_{ij}$, which satisfies $(xy)^* = y^*x^*$, and the coproduct $\Delta(u_{ij}) = \sum_{k=1}^n u_{ik} \otimes u_{kj}$. For details, see the exact construction in \cite{corey2023quantum}. 
The most important fact for us is that the defining ideal $I_n$ is closed under the involution. 
Thus it can be ommited for the ideal membership problem on $I_n$.
\begin{lemma}
    The generating set of the ideal $I_n$ is self-adjoint, i.e.
    \[
      \{ \rs_i^*,\cs_i^*,\inj_{ijk}^*,\wel_{ijk}^*,\ip_{ij}^* \} = \{ \rs_i,\cs_i,\inj_{ijk},\wel_{ijk},\ip_{ij} \}
    \]
    resulting in $I_n$ being closed under the $^*$-involution.
\end{lemma}
\begin{proof}
  Both $\rs_i$ and $\cs_i$ are self-adjoint, since they are sums of self-adjoint elements. Similarly,
   \begin{align*}
    &\ip_{ij}^* = (u_{ij}^2 - u_{ij})^* = u_{ij}^*u_{ij}^* - u_{ij}^* = u_{ij}u_{ij} - u_{ij} = \ip_{ji}\\
    &\wel_{ijk}^* = (u_{ij}u_{ik})^* = u_{ik}u_{ij} = \wel_{ikj}
    \end{align*}
    and 
    $\inj_{ijk}^* = (u_{ij}u_{kj})^* = u_{kj}u_{ij} = \inj_{kji}$.
\end{proof}
This opens up the possibility of solving the ideal membership problem, for example to get a commutativity result, by solving it in free associative algebra using Gröbner bases.

Before continuing with the theory of Gröbner bases, we will introduce another type of involution on $\C\langle n^2 \rangle$, which will help us in the proof of \Cref{thm:finite_gb}.

%\subsection{Transposition}

We investigate symmetries arising due to the matrix structure \eqref{eq:linear_preorder} and the specific shape of relations in \eqref{eq:def_quantumGroup}. 
For this define a second (multiplicative) involution. 
Note that this is not the $*$-involution  from the quantum group. 
    Let    \begin{align}\label{eq:def_symInvol}{\left(\cdot\right)}^{\times}:\C\langle n^2\rangle&\rightarrow\C\langle n^2\rangle\\
    u_{ij}&\mapsto u_{ji}\nonumber
    \end{align}
    be an algebra homomorphism, uniquely determined though the universal property. 

\begin{lemma}\label{lem:technicalDetails_SymInvolution}
  The map from \eqref{eq:def_symInvol} has the following properties.
    \begin{multicols}{2}
    \begin{enumerate}[i)]
        \item\label{lem:technicalDetails_SymInvolution1} The map \eqref{eq:def_symInvol} is an involution. 
        \item\label{lem:technicalDetails_SymInvolution2} If $v$ divides $w$ then $v^\times$ divides $w^\times$. 
        % \item $f\xrightarrow{G}r$ and $g^\times \in G$ for all $g\in G$, then $f^\times\xrightarrow{G}r^\times$. 
        \item\label{lem:technicalDetails_SymInvolution3} $\rs_j^\times=\cs_j$ 
        %\item $\cs_j^\times=\rs_j$  
        \item $\ip_{ij}^\times=\ip_{ji}$  
        \item $\wel_{ijk}^\times=\inj_{jik}$  
        %\item $\inj_{ijk}^\times=\wel_{jik}$
        \item $\rinj_{kj}^{\times}=\rwel_{kj}$
    %\item $\rwel_{kj}^{\times}=\rinj_{kj}$
    \end{enumerate}
    \end{multicols}
\end{lemma}
\begin{proof}
    Part \ref{lem:technicalDetails_SymInvolution1} follows with ${(u_{i,j}^\times)}^\times=u_{j,i}^\times=u_{i,j}$. For part \ref{lem:technicalDetails_SymInvolution2} let $w=avb$, then $w^{\times}={(avb)}^{\times}={a}^{\times}{v}^{\times}{b}^{\times}$. For part \ref{lem:technicalDetails_SymInvolution3}, 
    $$\rs_j^\times={\left(\sum_{1\leq \alpha\leq n}u_{j\alpha}-1\right)}^\times=\sum_{1\leq \alpha\leq n}u_{j\alpha}^\times-1=\sum_{1\leq \alpha\leq n}u_{\alpha j}-1=\cs_j.$$ %and therefore in part \ref{lem:technicalDetails_SymInvolution8}, $$(\lm(\rs_j))^\times=u_{j1}^\times=u_{1j}=\lm(\cs_j)=\lm(\rs_j^\times).$$ 
    The remaining parts are similar. 
\end{proof}

\section{Primer on non-commutative Gröbner bases}\label{sec:primGBs}

We recall the notion of non-commutative Gr\"obner bases. 
For a detailed introduction see the standard references \cite{bergman1978diamond}  or  \cite{mora1988groebner}. 
A more modern reference would be \cite{nordbeck2001canonical} with notation similar to ours.
For the entire section, we fix our free associative algebra 
$$R:=\mathbb{C}\langle n^2\rangle=\mathbb{C}\langle u_{ij}\mid 1\leq i,j\leq n\rangle$$
on the $n^2$ symbols $u_{ij}$ where $n \in \N, n \geq 4$. Furthermore, we choose the degree lexicographic order $<$ on $R$ with \DEF{row-wise linear preorder}, i.e., \begin{equation}\label{eq:linear_preorder}u_{i,j}>u_{i,j+1}\quad\text{ and }\quad u_{i,n}>u_{i+1,1}.
\end{equation}

\begin{example}We have $u_{11}>u_{12}>u_{21}>u_{22}$
and  $u_{12}u_{21}>u_{21}u_{12}>u_{12}$.
\end{example}

For any nonzero $f\in R$ let $\lm(f),\lc(f)$ and $\lt(f)$ denote its \DEF{leading monomial}, \DEF{leading coefficient} and \DEF{leading term} of $f$, respectively. 
 
\begin{definition}\label{def:gb_noncomm}
   A set $G\in R$ is a \DEF{(non-commutative) Gr\"obner basis} of a two-sided ideal $I\subseteq R$ if for every $f\in I\setminus\{0\}$ there exists $g\in G$ such that $\lm(g)$ divides $\lm(f)$.
\end{definition} We call any subset $G\subset R$  a Gr\"obner basis if it is a Gr\"obner basis of $\langle G\rangle$. 

\begin{remark}\label{rem:other_orderings_and_rewriting}
  An analogous version of Robbiano’s characterization \cite{Robbiano85} as in the commutative setting is not possible for free algebras, e.g. \cite{green1998non,hermiller1999monomial}. Note that the lexicographic order is not monomial. For this reason, we investigate its extension \eqref{eq:linear_preorder}. In particular, if a set is a Gr\"obner basis (\Cref{def:gb_noncomm}) then it is always with respect to our fixed order.
\end{remark}
\begin{example}\label{ex:injwel}
   The set of orthogonal relations $G:=
\{\inj_{ikj},\wel_{kij}\mid1\leq i,j,k\leq n,i\not=j\}$ is a finite  Gr\"obner basis since it is generated only by monomials. 
\end{example}

For any monomials $w,v\in R$ we say that $v$ \DEF{divides} $w$ if there are monomials $a,b\in R$ such that $w=avb$. 
With this we can formulate the \DEF{reduction algorithm}\footnote{The recursive call on the \emph{lower order terms} in \Cref{alg:reduction} is sometimes called \emph{tail reduction} and results in unique Gr\"obner bases in the sense of \Cref{prop:buchberger_correct}.} modulo a set of generators $G\subseteq R$, \Cref{alg:reduction}. 

\begin{algorithm}
\caption{Reduction algorithm}\label{alg:reduction}
\begin{algorithmic}
\Procedure{$\mathsf{NF}$}{$r,G$}
\While{$\exists g\in G\setminus\{0\}:\lm(g)\text{ divides }\lm(r)$}
\State choose $a,b$ with $a\lm(g)b=\lm(r)$
\State $r \gets r-\frac{\lc(r)}{\lc(g)}agb$
\EndWhile\label{euclidendwhile}
\State \textbf{return} $\lt(r)+\mathsf{NF}(r-\lt(r),G)$
\Comment{recursive tail reduction}
\EndProcedure
\end{algorithmic}
\end{algorithm}

Note that its output $\mathsf{NF}(r,G)$ is not unique in general since the algorithm requires choices. However, if $G$ is a Gr\"obner basis then the following result, also known as \emph{Bergman's diamond lemma}, shows that for every $r\in R$ the output $\mathsf{NF}(r,G)$ is uniquely determined. In this case  we call the output of the reduction algorithm  \DEF{normal form} of $r$ modulo $G$.

\begin{theorem}[Bergman \cite{bergman1978diamond}]\label{thm:bergman_diamond}
For any subset $G\subseteq R$, the following statements are equivalent. 
\begin{enumerate}[i)]
    \item $G$ is a Gr\"obner basis. 
    \item The output of the reduction algorithm $\mathsf{NF}(f,G)$ is unique for every $f\in R$. 
    \item The set of reduced monomials
    $$\{\mathsf{NF}(w,G)\mid w\text{ monomial }\}$$
     is a $\C$-basis of the factor algebra $\faktor{R}{\langle G\rangle}$ when considered as a vector space. 
\end{enumerate}
\end{theorem}

As an imitiate consequence of     \Cref{thm:bergman_diamond} we obtain the following. 

\begin{corollary}\label{cor:decision_if_finite_GB}
If $G$ is a finite Gr\"obner basis then the word problem in $\faktor{R}{\langle G\rangle}$ is decidable. 
\end{corollary}

Buchberger's algorithm is the essential tool for computing for a given input $F\subseteq R$ a Gr\"obner basis of the two-sided ideal $\langle F\rangle$. 
Note that it has no termination guarantees. However, if it terminates, then it provides a finite Gr\"obner basis and thus a decision procedure for the word problem in the associated factor algebra, \Cref{cor:decision_if_finite_GB}. 

We recall the essential criterion for Buchberger's algorithm. 
For every non-zero $f\in R$ and $G\subseteq R$ we call $(v,g,w)\in R^{3\ell}$ a \DEF{Gr\"obner representation}, if 
\begin{equation}\label{eq:GB_rep}
    f =\sum_{1\leq i\leq \ell}v_ig_iw_i
\end{equation}
with $\lm(f)\geq \lm(v_ig_iw_i)$ and $g_i\in G$ for all $1\leq i\leq \ell$. 
Similarly as for S-polynomials in the commutative Gr\"obner basis theory we have to enlarge our generating set so that  \eqref{eq:GB_rep} is always true. 
Two monomials $w,v\in R$ have an \DEF{overlap} if there are monomials $a,b\in R$ such that one of the conditions
\begin{align}
&wa=bv\quad\text{ or }\quad aw=vb\label{def:overlapOfWords2}%,\\
 %&w=aw'b,\text{ or 
%}\label{def:overlapOfWords3}\\
% \intertext{ or }
   %&awb=w'\label{def:overlapOfWords4}
\end{align}
with $0<\len(a)\leq\len(v)$ and $0<\len(b)\leq\len(w)$
is satisfied.
For all non-zero $f,g\in R$ and monomials $a,b\in R$  we obtain the following \DEF{overlap relations} 
\begin{equation}\label{eq:overlap_relations}
\begin{cases}\frac{1}{\lc(f)}fa-\frac{1}{\lc(g)}bg&\text{if }\lm(f)a=b\lm(g)\\
\frac{1}{\lc(f)}af-\frac{1}{\lc(g)}gb&\text{if }a\lm(f)=\lm(g)b
\end{cases}
 \end{equation}
with $a,b$ according to \eqref{def:overlapOfWords2}. Let $\mathcal{O}(f,g)$ denote the set of all those overlap relations. Note that this set is always finite. Similarly, if $\lm(f)$ divides $\lm(g)$ then we call   
\begin{align*}
&\frac{1}{\lc(f)}afb-\frac{1}{\lc(g)}g
    \end{align*}
 with $a\lm(f)b=\lm(g)$ a  \DEF{division relation}.  
With this we obtain the following non-commutative version of  \DEF{Buchberger's criterion} with finitely many conditions.
\begin{proposition}\label{prop:GB_via_GBrep}A subset $G\subset R$ is a Gr\"obner basis if and only if each overlap and division relation of any $f,g\in G$ has a Gr\"obner representation in $G$.
\end{proposition}
 We refer to  \cite[Proposition 6]{nordbeck2001canonical} for a proof of \Cref{prop:GB_via_GBrep}.

\begin{example}\label{ex_ip_is_groebner}
   The set of all idempotent relations $G:=\{\ip_{ij}\mid1\leq i,j\leq n\}$ is a Gr\"obner basis. We have no non-zero division relations and the only possible overlap relations in $G$ are self-overlaps $u_{ij}\ip_{ij}-\ip_{ij}u_{ij}=0$ according to \eqref{def:overlapOfWords2}, which have a trivial Gr\"obner  representation \eqref{eq:GB_rep}.  
\end{example}

We call a set $F\subseteq R$ 
\DEF{reduced} if none of the $\lm(f)$ divides $\lm(g)$ for $f,g\in F$ with $f\not=g$. 
We call  $F$ 
\DEF{tail-reduced} if $\mathsf{NF}(f,F\setminus\{f\})=f$ for all $f\in F$. Clearly, a tail-reduced set is also reduced. With a simple interreduction, \Cref{alg:interreduce}, we obtain for every input set $F$ a tail reduced set, denoted by $\mathsf{interreduce}(F)$.

\begin{algorithm}
\caption{Interreduction}\label{alg:interreduce}
\begin{algorithmic}
\Procedure{$\mathsf{interreduce}$}{$F$}

\While{$\exists f\in F$ with $r:=\mathsf{NF}(f,F\setminus\{f\})\not=f$}
\State replace $f$ by $r$ in $F$
\EndWhile\label{euclidendwhile2}
\State \textbf{return} $F\setminus\{0\}$
\EndProcedure
\end{algorithmic}
\end{algorithm}

We call a set $F\in R$ \DEF{monic} if $\lm(f)=1$ for all $f\in F$. 

\begin{proposition}\label{prop:uniqueniess_of_gbs}
    Every ideal has a unique Gröbner basis that is tail-reduced and monic. 
\end{proposition}

 We refer to  \cite[Proposition 14]{nordbeck2001canonical} for its proof. Finally, we can recall Buchberger's algorithm from \cite{mora1988groebner}, \Cref{alg:cap}. If it terminates, then its output is a Gr\"obner basis, \Cref{prop:buchberger_correct}. For a proof, compare \cite[Proposition 13]{nordbeck2001canonical}.

\begin{algorithm}
\caption{Buchberger's algorithm}\label{alg:cap}
\begin{algorithmic}
\Procedure{$\mathsf{Buchberger}$}{$F$}
\State $F \gets \mathsf{interreduce}(F)$
\While {$\exists f,g\in F\;\exists h\in\mathcal{O}(f,g)$ with $0\not= \mathsf{NF}(h,F)$}
\State $F \gets \mathsf{interreduce}(F\cup\mathcal{O}(f,g))$
\EndWhile
\State \textbf{return} $F$
\EndProcedure
\end{algorithmic}
\end{algorithm}

\begin{proposition}\label{prop:buchberger_correct}
    Let $F\subseteq R$ be finite set such that the two-sided ideal $\langle F\rangle$ has a finite  Gr\"obner basis. Then Buchberger's algorithm terminates and provides the monic and tail-reduced Gr\"obner basis of $\langle F\rangle$. 
\end{proposition}

The following result is already an important part of the proof of our main result, \Cref{thm:finite_gb}. Its proof illustrates interreduction (\Cref{alg:interreduce}) and the main proof strategy using Gr\"obner certificates  constructed via a  suitable chain of reductions (\Cref{rem:reduction4GroebnerCertificates}) that are valid for all sizes $n$. 

\begin{proposition}\label{ex:linRelations}
 The set $G:=\{\rs_2,\dots,\rs_n,\cs_1,\dots,\cs_n\}$ is the a reduced and  monic Gr\"obner basis of the ideal generated by all row and column sums, $\langle\rs_1 \cup G\rangle$.  
 \end{proposition}
 \begin{proof}
 For this we see that  $\lm(\rs_i)=u_{i1}$ and $\lm(\cs_i)=u_{1i}$ for every $i$,  hence $\len(\lm(g))=1$ for every $g\in G$. Therefore there are no overlaps in $G$. With $\lm(f)\not=\lm(g)$ for all $f,g\in G$,  we see that $G$ is reduced.  With 
 \begin{equation}\label{eq:reduction_rs_1}\rs_1=\sum_{1\leq i\leq n}\cs_i-\sum_{2\leq j\leq n}\rs_j\in\Span(G)
 \end{equation}
we obtain a Gr\"obner representation \eqref{eq:GB_rep} since $\lm(\rs_1)=u_{11}\geq u_{ij}$ for all $1\leq i,j\leq n$. 
 Note that $G$ is not tail-reduced. For instance,  $\lm(\rs_2)$ divides $\lm(\cs_1-\lt(\cs_1))$. 
\end{proof}
While \Cref{eq:reduction_rs_1} is easy to check by hand, it already foreshadows the computational requirements for the less trivial reductions later in this text.

\begin{remark}\label{rem:reduction4GroebnerCertificates}
We use the classical notation $f\xrightarrow{F}r$ when  $f$ reduces to $r$ modulo a set of divisors $F$. The Gr\"obner representation in \eqref{eq:reduction_rs_1} for instance is constructed via the reduction
$$\rs_1\xrightarrow{\cs_1}\rs_1-\cs_1\xrightarrow{\cs_2}\dots\xrightarrow{\cs_n}\rs_1-\sum_{1 \leq i\leq n}\cs_i\xrightarrow{\rs_2}\dots\xrightarrow{\rs_n}0.$$
We also use the  short-hand notation  $\rs_1\xrightarrow{G}0$, where $G$ is from \Cref{ex:linRelations}. 
\end{remark}

\begin{example}\label{ex:constructionOfGBcertificate}
    We illustrate \Cref{ex:linRelations} and \Cref{rem:reduction4GroebnerCertificates} for $n=3$. Then, 
    \begin{align*}
        \rs_1=u_{11}+u_{12}+u_{13}+1
        &\xrightarrow{\cs_{1}}
        u_{11}+u_{12}+u_{13}-1
        -
        \cs_{1}
        =
        u_{12}+u_{13}
        -
        u_{21}-u_{31}\\
        &\xrightarrow{\cs_{2}}
       u_{13}
        -
        u_{21}-u_{31}-u_{22}-u_{23}+1\\
        &\xrightarrow{\cs_{3}}
        u_{21}-u_{31}-u_{22}-u_{23}-u_{32}-u_{33}+2\\
        &\xrightarrow{\rs_{2}}
        -u_{31}-u_{32}-u_{33}+1\\
        &\xrightarrow{\rs_{3}}0,
    \end{align*}
hence the Gr\"obner representation 
$\rs_1=\cs_1+\cs_2-\rs_1-\rs_2-\rs_1
$. We refer to \Cref{ex:row_sums_col_sums_formaliz} for a complete explanation of how we verify \eqref{eq:reduction_rs_1} independently of $n$.  
\end{example}

In the following situation, the union of two Gr\"obner bases is again a Gr\"obner bases.  

\begin{proposition}\label{lem:smallGSs} The set 
$G=\{\ip_{ij},\inj_{jik},\wel_{ijk}\mid j\not=k\}$ is a Gr\"obner basis. 
\end{proposition}
\begin{proof}
Overlaps of $\inj$ and $\ip$ are covered in  \Cref{ex:injwel,ex_ip_is_groebner}. Overlaps between $\ip$ and $\inj$ are of the form 
        \[
            \ip_{ij}u_{jk} - u_{ik}\inj_{ikj} = u_{ik}u_{jk} \xrightarrow{\inj_{ikj}} 0
        \]
for $1\leq i,j,k\leq n$ with $i\not=j$. The remaining overlaps are analogous. 
\end{proof}

However, in general, Gr\"obner bases are not closed under union. We show in \Cref{thm:finite_gb} that the union of Gr\"obner bases from \Cref{ex:linRelations,lem:smallGSs} is almost a Gr\"obner basis of $I_n$, when including the extra relations from \Cref{def:rinjRwel}. 
This is by no mean obvious and we devote the entire \Cref{sec:main_result} for its proof. 
\par
Essential for this are certain compatibility laws regarding transposition and leading momomials of our quantum symmetry relations.

\begin{corollary}
\label{lem:technicalDetails_SymInvolution8}$(\lm(f))^\times=\lm(f^\times)$ for $f\in\{\rs_j, \cs_j,\ip_{ij},\wel_{ijk},\inj_{ijk},\rinj_{kj},\rwel_{kj}\}$
\end{corollary}
\begin{proof}
    With \Cref{lem:technicalDetails_SymInvolution} we have $\lm(\rs_j^\times)=\lm(\cs_j)=u_{1j}=u_{j1}^\times=(\lm(\rs_1))^\times$ for all $j$. The remaining  relations can be treated analogously. 
\end{proof}

\begin{example}
  Note that transposition is not always compatible with the ordering in the sense of \Cref{lem:technicalDetails_SymInvolution} and \Cref{lem:technicalDetails_SymInvolution8}. For instance with   $f:=u_{12}+u_{31}$ we have $$(\lm(f))^\times=u_{12}^\times=u_{21}\not=u_{13}=\lm(u_{21}+u_{13})=\lm(f^\times).$$
\end{example}

We close this section with another consequence of Bergman's diamond lemma, that the order of reductions can be exchanged in an arbitrary way. 

\begin{lemma}[Gr\"obner bases via extended relations]\label{lem:helpRelations}
  Any subset $G\subseteq R$ is a Gr\"obner basis if and only $G\cup\{ugv\}$ is a Gr\"obner basis, where $g\in G$ and $u,v$ are monomials.
    Furthermore,  $G$ is a Gr\"obner basis if and only if $G\cup\{f+g\}$
    is a Gr\"obner basis where $f,g\in G$ with $\lm(f)>\lm(g)$.
\end{lemma} 
%\begin{proof}This is well-known, but we give a proof for completeness. Clearly, if $G$ is a Gr\"obner basis then any $G\cup F$ with $F\subseteq\langle G\rangle$ is a Gr\"obner basis. It thus suffices to show the backwards directions. The first claim follows with interreduction since $ugv$ reduced to zero for  all monomials $u,v$ and $g\in G$. For the second claim, 
%    let $f,g,h\in G$ with $\lm(f)>\lm(g)$. %Then 
%    %$\lm(u(f+g)w)=u\lm(f)v$ and 
%    %$$u(f+g)v\xrightarrow{g}ufg.$$
%    If $f$ overlaps with $h$ via $\lm(af)=\lm(gb)$  for suitable monomials $a$ and $b$, then $(f,h)$ and $(f+g,h)$ yield the same overlap constellations, and the corresponding relations can be reduced to the same element, 
%    $$\frac1{\lc(f)}a(f+g)-\frac1{\lc(h)}hb\xrightarrow{g}\frac1{\lc(f)}af-\frac1{\lc(h)}hb.$$
%\end{proof}

\section{Proof of \texorpdfstring{\Cref{thm:finite_gb}}{Theorem \ref{thm:finite_gb}}}
\label{sec:main_result}

%\begin{theorem}\label{thm:finite_gb}
%For all $n\geq 4$ the ideal $I_n$ has a finite, monic and reduced Gröbner basis 
 %explicitly provided in . 
%\end{theorem}
In this section we prove our main result, i.e., we show that $G_n$ is a Gr\"obner basis for the ideal of quantum symmetries $I_n$. 
The basic proof strategy is to apply Buchberger's algorithm in a slightly modified form and independently of the size  $n$. We organize the proof into the following subsections. 
\par 
In \Cref{subsec:proof_first_interreduce} we reduce our generators from $I_n$. This is a classical  interreduction (\Cref{alg:interreduce}) without the recursive tail-reduction. We denote the resulting reduced set by $F_n$. 
In \Cref{subsec:mainproof_overlaps_bg} we determine all overlap relations  in $F_n$. When including those, and computing its reduced set, we already result in $G_n$. Finally, in \Cref{subsec:proof_veryBadGuys}, we apply Buchberger's criterion (\Cref{prop:GB_via_GBrep}), i.e., we show that all overlap relations in $G_n$ have a Gr\"obner representation. Note that these representations have to be generalizable for arbitrary $n \geq 4$. %For this we extend our quantum relations according to \Cref{lem:helpRelations}. 

\subsection{An interreduced generating set for \texorpdfstring{$I_n$}{In}}
\label{subsec:proof_first_interreduce}

We start with the preliminary set of generators for the ideal $I_n$, 
\begin{equation}\label{eq:gens_of_In}F''_n:=\left\{\rs_i,\cs_i,\ip_{ij}, \inj_{ijk},\wel_{ijk}\;\begin{array}{|l}\;1\leq i,j,k\leq n\\\text{ with }j\not=k\end{array}\right\}
\end{equation}
as it is defined in the introduction, that is $I_n=\langle F''_n\rangle$. Furthermore, we recall the reduced orthogonal relations $\rinj_{kj}$ and $\rwel_{kj}$ and show that their names are justified in the following sense.  
\begin{lemma}\label{lem:rinj_rwel_name_just}
\begin{enumerate}[i)]
\item If $2\leq j,k\leq n$ and $j\not= k$ then  
\begin{align*}
    \rinj_{jk} &=\inj_{j1k} - \rs_{j} u_{k1} + u_{j2} \rs_{k} - \inj_{j2k},\\
    \rwel_{jk} &= \wel_{1jk} - \cs_{j} u_{1k} + u_{2j} \cs_{k} - \wel_{2jk}
\end{align*} 
and both  $\rwel_{jk}$ and $\rinj_{jk}$ are reduced modulo $F''_n$. 
\item For all other cases of $j$ and $k$, both $\wel_{jk}$ and $\inj_{jk}$ reduce to zero modulo $F''_n$.   
\end{enumerate}
\end{lemma}
\begin{proof}
  For $k=1$ and all $i,j\not=1$, 
    \begin{align*}
        \inj_{1ij}=u_{1i}u_{ji}\xrightarrow{\cs_j}-\sum_{\alpha \neq 1}u_{\alpha i}u_{ji}+u_{ji}\xrightarrow{\inj_{\alpha ij}}-\ip_{ji} \\
        \inj_{ji1}=u_{ji}u_{1i}\xrightarrow{\cs_j}-\sum_{\alpha \neq 1}u_{ji}u_{\alpha i}+u_{ji}\xrightarrow{\inj_{j\alpha i}}-\ip_{ji}.
    \end{align*}
For arbitrary $k$ and $k\not=j\not=1$, 
    \begin{align*}
        \inj_{k1j}= u_{k1}u_{j1}
        \xrightarrow{\rs_k}&u_{k1}u_{j1}-\left(\sum_{\alpha \in E}u_{k\alpha }-1\right)u_{j1} 
        =-\sum_{\alpha \not=1}u_{k\alpha}u_{j1} + u_{j1}\\
        \xrightarrow{\rs_j}& -\sum_{\alpha \not=1}u_{k\alpha}u_{j1} + u_{j1} - (-u_{k2})\left(\sum_{\alpha \in E}u_{j\alpha} -1\right)\\
        =& \sum_{\alpha \neq 1}\underline{u_{k2}u_{j\alpha}} -\sum_{\alpha \neq 1,2}u_{k\alpha}u_{j1}+ u_{j1} - u_{k2} \\
        \xrightarrow{\inj_{k2j}}& \sum_{\substack{\alpha \neq 1,2}}\underline{u_{k2}u_{j\alpha}} - \sum_{\alpha \neq 1,2}u_{k\alpha}u_{j1} + u_{j1} - u_{k2} = \mathsf{rinj}_{kj} 
\end{align*}
with leading monomial $\lm(\mathsf{rinj}_{kj})=u_{k2}u_{j3}$.  For $k \neq 1$ this is a reduced version of $\inj_{k1j}$. 
If $k=1$ we have
\begin{align}
  \inj_{11j}=u_{11}u_{j1}\xrightarrow{\rs_1}\sum_{s\not=1}-u_{s1}u_{j1}+u_{j1}=\sum_{s\neq1,s\not=j}\inj_{s1j}-\ip_{j1}. 
\end{align}
The analogous statements for $\wel$ follow from \Cref{lem:technicalDetails_SymInvolution}, e.g., 
$$\rwel_{jk}=\rinj_{jk}^\times=\inj_{j1k}^\times - \rs_{j}^\times u_{k1}^\times + u_{j2}^\times \rs_{k}^\times - \inj_{j2k}^\times=\wel_{1jk} - \cs_{j} u_{1k} + u_{2j} \cs_{k} - \wel_{2jk}.$$
\end{proof}

We reduce $F''_n$ and obtain the following smaller generating set $F'_n$. 

\begin{lemma}\label{lem:almost_reduced_generating_set_In}
  For all $n\geq 4$,
\[
  F_n':=\{\cs_1\}\cup \{ \cs_i,\rs_i,\ip_{ij},\inj_{jik},\wel_{ijk},\rinj_{jk},\rwel_{jk}\mid i,j,k\not=1 \text{ and } j\not=k\}
\]
is a generating set of $I_n$.
\end{lemma}
\begin{proof}
For every $j$, 
 \begin{align*}
        \ip_{1j}=u_{1j}u_{1j}-u_{1j}
                        &\xrightarrow{\rs_{j}} u_{ij}u_{ij}-u_{ij} - u_{ij}\left(\sum_{\alpha \in E}u_{i\alpha}-1\right)=-\sum_{\alpha \not=j}u_{ij}u_{i\alpha}\xrightarrow{\wel_{ij\alpha}}0,
    \end{align*}
    and similarly for $\ip_{j1}$ with \Cref{lem:technicalDetails_SymInvolution}. 
  We saw in  \Cref{ex:linRelations} how to reduce $\rs_1$ modulo $F_n''\setminus\{\rs_1$\}.   The remaining follows with \Cref{lem:rinj_rwel_name_just}. 
\end{proof}

\begin{proposition}\label{prop:red_gens_Fn}
 With $F_n := F_n' \setminus \{\rwel_{23}\}$ we obtain a reduced generating set of $I_n$. Its cardinality is given by the cubic polynomial $\#F_n=2n^3 - 5n^2 + 4n-1$. 
\end{proposition}

\begin{proof}
We have 
\begin{itemize}
  \item $2n-1$ row and column sum relations,
  \item $(n-1)^2 \cdot (n-2)$ relations of the form $\wel_{ijk}$ with $i,j,k\neq 1$ and $k\neq j$,
  \item $(n-1)^2 \cdot (n-2)$ relations of the form $\inj_{jik}$ with $i,j,k\neq 1$ and $k\neq j$,
\item $(n-1)(n-2)$ relations of the form $\mathsf{rwel}_{kj}$ with $k,j\neq 1$ and $k\neq j$,
  \item $(n-1)(n-2)$ relations of the form $\mathsf{rinj}_{kj}$ with $k,j\neq 1$ and $k\neq j$, and 
  \item $(n-1)^2$ relations of the form $\ip_{ij}$ with $i,j\neq 1$.
\end{itemize}
In total we obtain $2n-1 + 3(n-1)^2 + 2(n-2)(n-1)^2=2n^3 - 5n^2 + 4n$ and thus the claimed cardinality when omitting $\rwel_{23}$. Clearly $F'_n$ is a generating set with \Cref{lem:almost_reduced_generating_set_In}. We devote  
\Cref{subsection:reduction_rwel23} for a proof that 
$\rwel_{23}$ reduces to zero modulo $F_n$, hence also $F_n$ is a generating set of $I_n$. This remaining reduction of $\rwel_{23}$ is elementary but  quite long and  technical. Therefore, we  postpone it to the end of this manuscript together with the associated  computational machinery. 
\end{proof}

\begin{figure}
   \begin{align*}
      \bg{1}_{kji} &= \inj_{k2j} u_{i3}-u_{k2}\rinj_{ji} &
      \bg{2}_{kji} &=u_{k2} \inj_{j3i}-\rinj_{kj}u_{i3} \\
      \bg{3}_{kj}  &=\ip_{2k} u_{3j}-u_{2k}\rwel_{kj} &
      \bg{4}_{kj} &=u_{2k}\ip_{3j}-\rwel_{kj}u_{3j} \\
      \bg{5}_{kj} &=\ip_{k2}u_{j3}-u_{k2}\rinj_{kj} &
      \bg{6}_{kj} &=u_{k2}\ip_{3j}-\rinj_{kj}u_{3j}\\
      \bg{7}_{kji} &=\wel_{2kj}u_{3i}-u_{2k}\rwel_{ji} &
      \bg{8}_{kji} &=u_{2k}\wel_{3ji}-\rwel_{kj}u_{3i} \\
      \bg{9}_{kj} &=\rinj_{k2}u_{3j}-u_{k2}\rwel_{3j} &
      \bg{10}_{kj} &=u_{2k}\rinj_{3j}-\rwel_{k2}u_{j3} \\
      \bg{11}_{kji} &=\inj_{kj2}u_{3i}-u_{kj}\rwel_{ji} &
      \bg{12}_{kji} &=u_{2k}\inj_{3ji}-\rwel_{kj}u_{ij} \\
      \bg{13}_{kji} &=\wel_{kj2}u_{i3}-u_{kj}\rinj_{ki} &
      \bg{14}_{kji} &=u_{k2}\wel_{j3i}-\rwel_{kj}u_{ji} 
   \end{align*} 
\caption{List of all possible overlaps between families in $F_n$.}
  \label{def:G0_overlaps}
\end{figure}
\subsection{Construction of the Gr\"obner basis  \MathInHeader{G_n}}\label{subsec:mainproof_overlaps_bg}

\begin{figure}
  \centering
  \begin{tikzpicture}[scale=3, every node/.style={font=\sffamily\bfseries, inner sep=1pt}, every edge/.style={draw, thick}]

    % Die Knoten im Pentagon
    \node[draw, shape=rectangle, inner sep=5pt] (ip)  at (91:1)   {$\ip$};
    \node[draw, shape=rectangle, inner sep=5pt] (inj) at (162:1)  {$\inj$};
    \node[draw, shape=rectangle, inner sep=5pt] (wel) at (234:1)  {$\wel$};
    \node[draw, shape=rectangle, inner sep=5pt] (rwel) at (306:1) {$\rwel$};
    \node[draw, shape=rectangle, inner sep=5pt] (rinj) at (18:1) {$\rinj$};

    % Die Kanten zwischen den Knoten
    \path[-] (ip)  edge node[sloped, above, font=\small, draw=none] {\hyperref[lem:bg5bg6]{\textcolor{black}{$\bg{5},\bg{6}$}}} (rinj);
    \path[-] (ip)  edge node[sloped, above, font=\small, draw=none] {\Cref{lem:smallGSs}} (inj);
    \path[-] (ip)  edge node[sloped, above, font=\small, draw=none, pos=0.79] {Prop. \ref{lem:smallGSs}} (wel);
    \path[-] (ip)  edge node[sloped, above, font=\small, draw=none, pos=0.21] {\hyperref[lem:overlaps_bg34]{\textcolor{black}{$\bg{3},\bg{4}$}}} (rwel);

    \path[-] (inj) edge node[sloped, below, font=\small, draw=none] {\Cref{ex:injwel}} (wel);
    \path[-] (inj) edge node[sloped, above, font=\small, draw=none] {\hyperref[lem:bg11bg12]{\textcolor{black}{$\bg{11},$}}}
      node[sloped, below, font=\small, draw=none] {\hyperref[lem:bg11bg12]{\textcolor{black}{$\;\bg{12}$}}}
      (rwel);
    \path[-] (inj) edge node[sloped, above, font=\small, draw=none, pos=0.21] {\hyperref[lem:bg1bg2]{\textcolor{black}{$\bg{1},\bg{2}$}}} (rinj);

    \path[-] (wel) edge node[sloped, above, font=\small, draw=none] {\hyperref[lem:overlap78]{\textcolor{black}{$\bg{7},\bg{8}$}}} (rwel);
    \path[-] (wel) edge node[sloped, below, font=\small, draw=none, pos=0.8] {\hyperref[lem:bg13bg14]{\textcolor{black}{$\bg{14}$}}}
      node[sloped, above, font=\small, draw=none, pos=0.82] {\hyperref[lem:bg13bg14]{\textcolor{black}{$\bg{13},$}}}(rinj);

    \path[-] (rwel) edge node[sloped, below, font=\small, draw=none] {\hyperref[lem:overlaps_gb910]{\textcolor{black}{$\bg{9},\bg{10}$}}} (rinj);

  \end{tikzpicture}
  \caption{Graph of all possible overlap parings between families in $F_n$.}
  \label{fig:pentagonGraphG0}
\end{figure}
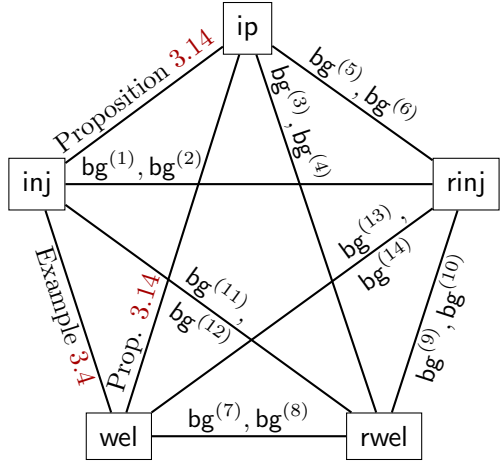

In this section we determine all overlap relations in $F_n$, listed for illustration in \Cref{def:G0_overlaps}. They are grouped as overlap relations of certain pairs of families in $F_n$. We illustrate the underlying parings of those families in a graph, \Cref{fig:pentagonGraphG0}.
Note that row and column sums cannot produce overlaps since their leading monomials are of length $1$.

Only the $2n(n-2)(n-3)-1$ overlap relations from two families, 
\begin{equation}\label{eq:reducedBadGuys_new}B_n:=\left\{\bg{s}_{kji}\; \begin{array}{|l}s\in\{2,8\},k,j,i\in[2,n]\\
\text{where }i\not=j\not=k,
\\
(k,j)\not=(2,3)\not=(j,i)\text{  and } \\(s,k,j,i)\not=(8,2,4,3)\end{array}\right\}\end{equation} do not reduce to zero after a suitable reduction. Furthermore, we show that each relation in the latter is already  reduced.
Note that $B_n$ is precisely the union of the second and third set from the  disjoint union of our Gr\"obner basis $G_n$ in \Cref{thm:finite_gb}, i.e. $G_n=F_n\cup B_n$ for all $n\geq 4$. 
We start with the first paring of $\inj$ and $\rinj$, illustrated also in \Cref{fig:overlaps_bg12}.

\begin{lemma}[$\inj$ and $\rinj$]
  \label{lem:bg1bg2}
 There are two types of overlaps for $\inj$ and $\rinj$,
  \begin{align*}
      \bg{1}_{ijk}& := \inj_{i2j} u_{k3}-u_{i2}\rinj_{jk}\quad \text{ for all } i,j,k\in[2,n] \text{ with } i \neq j \wedge j \neq k,\text{ and}\\
      \bg{2}_{ijk}& :=u_{i2} \inj_{j3k}-\rinj_{ij}u_{k3} \quad\text{ for all } i,j,k\in[2,n] \text{ with } j \neq k \wedge i \neq j.
 \end{align*} 
 Every overlap relation of type $\bg{1}_{ijk}$ reduces to zero modulo $F_n$. Overlaps of type $\bg{2}_{ijk}$ are reduced modulo $F_n$ for all 
 $(i,j)\not=(2,3)$ or $(j,k)\not=(2,3)$.
\end{lemma}

\begin{proof}
  We start with a proof for $\bg{1} \xrightarrow{F_n} 0$. If $i,j,k\in[2,n]$ with $i \neq j \wedge j \neq k$, then 
      \begin{align*}
            \bg{1}_{ijk} = & - \sum_{\alpha > \hlmath{3}} u_{i2}u_{j2}u_{k\alpha} + \sum_{\alpha > 2} u_{i2}u_{j \alpha} u_{k1} - u_{i2}u_{k1} + u_{i2}u_{j2} \\
          \xrightarrow{\inj_{i2j}}&   \sum_{\alpha > 2} u_{i2}u_{j \alpha} u_{k1} - u_{i2}u_{k1} \\
          \xrightarrow{\rinj_{ij}u_{k1}}& \sum_{\alpha > 2} u_{i2}u_{j \alpha} u_{k1} - u_{i2}u_{k1} - \sum_{\alpha > 2} u_{i2}u_{j\alpha}u_{k1} + \sum_{\alpha > 2} u_{i\alpha}u_{j1}u_{k1} - u_{j1}u_{k1} + u_{i2}u_{k1} \\
          =&  \sum_{\alpha > 2}u_{i2}u_{j1}u_{k1} -u_{j1}u_{k1}
          \xrightarrow{\inj_{j1k}} 0.
          \end{align*}
For the second claim we observe that $\lm(\bg{2}_{ijk})=u_{i2}u_{j4}u_{k3}$ is not divisible by any monomial in $\lm(F_n)$. 
\end{proof}
\begin{figure}[htbp]
  \centering
  \begin{subfigure}[b]{0.45\textwidth}
    \centering
\begin{tabular}{cc *{2}{@{$\mskip\thickmuskip$}c@{$\mskip\thickmuskip$}c}}
  $\lm(\inj_{i2j}u_{k3}):$ & $u_{i2}$ & $\cdot$ & $u_{j2}$ &\vline & $u_{k3}$ \\ \cline{2-6}\cline{2-6}
  $\lm(u_{i2}\rinj_{jk}):$ & $u_{i2}$ &\vline & $u_{j2}$ & $\cdot$ & $u_{k3}$ \\
\end{tabular}
    \caption*{$\bg{1}$}
  \end{subfigure}
  \begin{subfigure}[b]{0.45\textwidth}
    \centering
\begin{tabular}{cc *{2}{@{$\mskip\thickmuskip$}c@{$\mskip\thickmuskip$}c}}
      $\lm(u_{i2}\inj_{j3k}):$ & $u_{i2}$ & \vline & $u_{j2}$ &$\cdot$& $u_{k3}$ \\ \cline{2-6}\cline{2-6}
    $\lm(\rinj_{ij}u_{k3}):$ & $u_{i2}$ &$\cdot$ & $u_{j3}$ & \vline & $u_{k3}$ \\
    \end{tabular}
    \caption*{$\bg{2}$}
  \end{subfigure}
\caption{The two  possible overlaps of the paring $\inj$ and $\rinj$, resulting in the two families of overlap relations, $\bg{1}$ and $\bg{2}$.}
\label{fig:overlaps_bg12}
\end{figure}
 
We cover the two remaining cases $(i,j)=(2,3)$ and $(j,k)=(2,3)$  that are not addressed in \Cref{lem:bg1bg2} later in the proof of \Cref{lem:gb2_243_minus_bg8_243}. For the moment it suffices to note that they are not included in $B_n$. We move on to the second pairing.

\begin{lemma}[$\ip$ and $\rwel$]\label{lem:overlaps_bg34}
There are two types of overlaps for $\ip$ and $\rwel$, 
\begin{align*}
      \bg{3}_{kj} &:=\ip_{2k} u_{3j}-u_{2k}\rwel_{kj} \quad \text{ for all } k,j\in[2,n] \text{ with } k \neq 2 \wedge k \neq j,\text{ and}\\
      \bg{4}_{kj}& :=u_{2k}\ip_{3j}-\rwel_{kj}u_{3j}\quad \text{ for all } k,j\in[2,n] \text{ with } j \neq 3 \wedge k \neq j.
\end{align*}
Both reduce to zero modulo $F_n$. 
\end{lemma}
\begin{proof}
  Set $k,j\in[2,n]$ with $k \neq 2 \wedge k \neq j$. Then, 
  \begin{align*}
    \bg{3}_{kj} &= \ip_{2k} u_{3j}-u_{2k}\rwel_{kj} \\
                &= - \sum_{\alpha > \hlmath{3}} u_{2k}u_{2k}u_{\alpha j} + \sum_{\alpha > 2} u_{2k}u_{\alpha k} u_{1j} - u_{2k}u_{1j} + u_{2k}u_{2k} - \hlmath{u_{2k}u_{3j}}\\
                &\xrightarrow{\inj_{2k\alpha}u_{1j}} - \sum_{\alpha > 3} u_{2k}u_{\alpha j} - u_{2k}u_{1j} + u_{2k}u_{2k} - u_{2k}u_{3j}\\
                &\xrightarrow{\ip_{2k},u_{2k}\cs_j} -\sum_{\alpha > 3} u_{2k}u_{\alpha j} + \sum_{\alpha \neq 2} u_{2k}u_{\alpha j} - u_{2k} + u_{2k} - u_{2k}u_{3j} = 0. %\\ 
  \end{align*}
 % $\bg{4}_{kj}$ is reduced in a 
  Similarly let  $k,j\in[2,n]$ with $j \neq 3 \wedge k \neq j$, then 
   \begin{align*}
            \bg{4}_{kj} %&= u_{2k}\ip_{3j}-\rwel_{kj}u_{3j} \\
            &= \sum_{\alpha > \hlmath{3}} u_{2k}u_{\alpha j}u_{3j} + \sum_{\alpha > 2} u_{\alpha k}u_{1j} u_{3j} - u_{1j}u_{3j} + u_{2k}u_{3j} - \hlmath{u_{2k}u_{3j}}
            \xrightarrow{\inj_{\alpha j 3 }, \inj_{1j3}} 0.
    \end{align*}
   % \incrit{Flesh this out some more, go into more Detail}.
\end{proof}

\begin{lemma}[$\ip$ and $\rinj$]\label{lem:bg5bg6}
 There are two types of overlaps for $\ip$ and $\rinj$, 
 \begin{align*}
      \bg{5}_{kj}&:=\ip_{k2}u_{j3}-u_{k2}\rinj_{kj} \quad\text{ for all } k,j\in[2,n] \text{ with } k \neq 2 \wedge k \neq j,\\
      \bg{6}_{kj}&:=u_{k2}\ip_{3j}-\rinj_{kj}u_{3j} \quad\text{ for all } k,j\in[2,n] \text{ with } j \neq 3 \wedge k \neq j.
 \end{align*}
Both reduce to zero modulo $F_n$.
 \end{lemma}
 \begin{proof}
     We can prove this directly via a suitable reduction. Alternatively, as in \Cref{lem:rinj_rwel_name_just}, we can use  \Cref{lem:technicalDetails_SymInvolution} for 
     $$\left(\bg{5}_{ijk}\right)^\times=\bg{3}_{ijk}\quad\text{ and }
     \quad\left(\bg{6}_{ijk}\right)^\times=\bg{4}_{ijk}$$ 
     and transpose the Gr\"obner representation in  
     \Cref{lem:overlaps_bg34}, which is again a Gr\"obner representation due to \Cref{lem:technicalDetails_SymInvolution8}. 
 \end{proof}

\begin{lemma}[$\wel$ and $\rwel$]\label{lem:overlap78}
 There are two types of overlaps for $\wel$ and $\rwel$,
 \begin{align*}
      \bg{7}_{ijk}&:=\wel_{2ij}u_{3k}-u_{2i}\rwel_{jk} \quad\text{ for all } i,j,k\in[2,n] \text{ with } i \neq 2 \wedge j \neq k,\text{ and}\\
      \bg{8}_{ijk}&:=u_{2i}\wel_{3jk}-\rwel_{ij}u_{3k}\quad \text{ for all } i,j,k\in[2,n] \text{ with } j \neq k \wedge i \neq j.
 \end{align*}
 Every overlap relation of type $\bg{7}_{ijk}$ reduces to zero modulo $F_n$. Overlaps of type $\bg{8}_{ijk}$ are reduced modulo $F_n$ for all 
 $(i,j)\not=(2,3)$ or $(j,k)\not=(2,3)$.
\end{lemma}
\begin{proof}
    This is again, using 
    $$\left(\bg{1}_{ijk}\right)^\times=\bg{7}_{ijk}\quad\text{ and }
     \quad\left(\bg{2}_{ijk}\right)^\times=\bg{8}_{ijk},$$ 
     the transposed version of \Cref{lem:bg1bg2}. 
\end{proof}

\begin{lemma}[$\rinj$ and $\rwel$]\label{lem:overlaps_gb910}
There are two types of overlaps for $\rinj$ and $\rwel$, 
 \begin{align*}
      \bg{9}_{kj}&:=\rinj_{k2}u_{3j}-u_{k2}\rwel_{3j}\quad \text{ for all } k,j\in[2,n] \text{ with } k \neq 2 \wedge j \neq 3,\text{ and}\\
      \bg{10}_{kj} &:=u_{2k}\rinj_{3j}-\rwel_{k2}u_{j3} \quad\text{ for all } k,j\in[2,n] \text{with } j \neq 3 \wedge k \neq 2.
 \end{align*} 
Both reduce to zero modulo $F_n$.
 \end{lemma}
\begin{proof}
We give the Gr\"obner representation and proof of its validity for $\bg{9}_{jk}$ in the provided git repository.
As seen above we obtain with transposition a Gr\"obner representation also for  $\bg{10}_{st}=-\left(\bg{9}_{st}\right)^\times$. 
\end{proof}

\begin{lemma}[$\inj$ and $\rwel$]\label{lem:bg11bg12}
  There are two types of overlaps for  $\inj$ and $\rwel$, 
 \begin{align*}
      \bg{11}_{kji} &:=\inj_{kj2}u_{3i}-u_{kj}\rwel_{ji}\quad\text{ for all } k,j,i \in [2,n] \text{ with } j\neq 2 \wedge j \neq i,\text{ and}\\
      \bg{12}_{kji}&:=u_{2k}\inj_{3ji}-\rwel_{kj}u_{ij} \quad\text{ for all } k,j,i \in [2,n] \text{ with } j \neq i \wedge k \neq j.
 \end{align*} 
 Both reduce to zero modulo $F_n$.
\end{lemma}
\begin{proof}
  For all  $k,j,i \in [2,n]$ with $j\neq 2 \wedge j \neq i$, 
  \begin{align*}
    \bg{11}_{kij}&=\inj_{ki2}u_{3j}-u_{ki}\rwel_{ij}\\&=-\sum_{\alpha>3}u_{ki}u_{2i}u_{\alpha j}+\sum_{\alpha>2}u_{ki}u_{\alpha i}u_{1j}-u_{ki}u_{1j}+u_{ki}u_{2i}&\\&\xrightarrow{\inj}\ip_{ki}u_{1j}+\inj_{ki2}. 
  \end{align*}
  Let $k,j,i \in [2,n]$ with $j \neq i \wedge k \neq j$, 
  \begin{align*}
    \bg{12}_{kij}&=u_{2k}\inj_{3ij}-\rwel_{ki}u_{ji}\\
                 &=\sum_{\alpha>3}u_{2k}u_{\alpha i}u_{ji}+\sum_{\alpha>2}u_{\alpha k}u_{1i}u_{ji}+u_{1i}u_{ji}-u_{2k}u_{ji}\\
                 &\xrightarrow{\wel}\begin{cases}
                   u_{2k}u_{ji}u_{ji}-u_{2k}u_{ji}&\text{if }j\not=2\\
                   u_{2k}u_{2i}&\text{if }j=2
                 \end{cases} \xrightarrow{\ip,\wel} 0. 
  \end{align*}
\end{proof}

\begin{lemma}[$\wel$ and $\rinj$]\label{lem:bg13bg14}
  There are two types of overlaps for $\wel$ and $\rinj$, 
 \begin{align*}
      \bg{13}_{kji}&:=\wel_{kj2}u_{i3}-u_{kj}\rinj_{ki}\quad \text{ for all } k,j,i \in [2,n] \text{ with } j\neq 2 \wedge k \neq i,\text{ and}\\
      \bg{14}_{kji} &:=u_{k2}\wel_{j3i}-\rwel_{kj}u_{ji} \quad\text{ for all } k,j,i \in [2,n] \text{ with  } i\neq 3 \wedge k \neq j.
 \end{align*} 
 Both reduce to zero modulo $F_n$.
  \end{lemma}
\begin{proof}
    This is again, using 
    $$\left(\bg{11}_{ijk}\right)^\times=\bg{13}_{ikj}\quad\text{ and }
     \quad\left(\bg{12}_{ikj}\right)^\times=\bg{14}_{ikj},$$ 
     the transposed version of \Cref{lem:bg11bg12}. 
\end{proof}

In order to determine all overlaps in $F_n$ it remains to consider the pairings of each class with itself.

\begin{lemma}[Self-overlaps]\label{lem:selfoverlaps_F1}
  There are no types in $F_n$  with overlaps that do not reduce to zero modulo $F_n$. 
\end{lemma}
\begin{proof}
    We refer to \Cref{ex:injwel,ex_ip_is_groebner}. The remaining cases are analogous. 
\end{proof}

\begin{proposition}\label{lem:gb2_243_minus_bg8_243}
\begin{enumerate}[i)]
    \item 
The union $G_n=F_n\cup B_n$ is a reduced set of generators for $I_n$.
\item Every overlap in $F_n$ is either contained in $B_n$ or reduces to zero modulo $G_n$. 
\end{enumerate}
\end{proposition}

\begin{proof}
The first part follows from \Cref{prop:red_gens_Fn}.
In Lemmas \ref{lem:bg1bg2} to \ref{lem:selfoverlaps_F1} we cover all overlaps in $F_n$, recorded in the graph of all parings, \Cref{fig:pentagonGraphG0}. 
\par
Once again, the reduction $\bg{8}_{243}$ and $\bg{2}_{23i}$ modulo $G_n$ is provided in the accompanying git repository. The latter covers the omitted reduction that is not provided by \Cref{lem:bg1bg2}.  
\par
The certificate for $\bg{2}_{k23}$ is analogous. With transposition we obtain a certificate for $$\left(\bg{2}_{23i}\right)^\times=\bg{8}_{23i}\quad\text{ and }\quad\left(\bg{2}_{k23}\right)^\times=\bg{8}_{k23},$$ 
and hence the omitted reductions not provided by  \Cref{lem:overlap78}. 
\end{proof}

\subsection{Remaining overlap relations in \MathInHeader{G_n}}\label{subsec:proof_veryBadGuys}

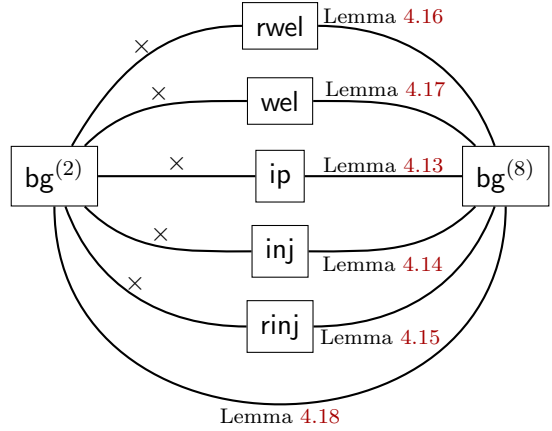
\begin{figure}
  \centering
  \begin{tikzpicture}[scale=3, every node/.style={font=\sffamily\bfseries, inner sep=1pt}, every edge/.style={draw, thick}]

    \node[draw, shape=rectangle, inner sep=5pt]   (ip)   at (0:0)   {$\ip$};
    \node[draw, shape=rectangle, inner sep=5pt, above of=ip] (wel) {$\wel$};
    \node[draw, shape=rectangle, inner sep=5pt, below of=ip] (inj) {$\inj$};
    \node[draw, shape=rectangle, inner sep=5pt, below of=inj]  (rinj)  {$\rinj$};
    \node[draw, shape=rectangle, inner sep=5pt, above of=wel]  (rwel)  {$\rwel$};
    \node[draw, shape=rectangle, inner sep=5pt, left of=ip, xshift=-2cm]  (bg2)  {$\bg{2}$};
    \node[draw, shape=rectangle, inner sep=5pt, right of=ip, xshift=2cm]  (bg8)  {$\bg{8}$};
    
    % Die Kanten zwischen den Knoten
\path[-, out=60, in=180] (bg2) edge node[ above, font=\small, draw=none] {$\times$} (rwel);
    \path[-, out=45, in=180] (bg2) edge node[ above, font=\small, draw=none] {$\times$} (wel);
    \path[-] (bg2) edge node[ above, font=\small, draw=none] {$\times$} (ip);
    \path[-, out=315, in =180] (bg2) edge node[ above, font=\small, draw=none] {$\times$} (inj);
    \path[-, out=290,in=180] (bg2) edge node[ above, font=\small, draw=none] {$\times$} (rinj);

    \path[-, out=270,in=270,looseness=1.5] (bg2) edge node[ below, font=\scriptsize, draw=none] {\Cref{lem:overlapsG0}} (bg8);

    \path[-, out=110,in=0] (bg8) edge node[pos=0.72, above, font=\scriptsize, draw=none,yshift=1.2mm] {\Cref{lem:bg8andrwel}} (rwel);
    \path[-, out=135,in=0] (bg8) edge node[pos=0.6, above, font=\scriptsize, draw=none,yshift=0.5mm] {\Cref{lem:bg8andwel}} (wel);

    \path[-] (bg8) edge node[sloped, above, font=\scriptsize, draw=none] {\Cref{lem:overlaps_bg8_ip}} (ip);
    \path[-, out=225, in=0] (bg8) edge node[pos=0.6, below, font=\scriptsize, draw=none, yshift=-0.5mm] {\Cref{lem:overlaps_bg8_inj}} (inj);
    \path[-, out=250, in=0] (bg8) edge node[pos=0.72, below, font=\scriptsize, draw=none,yshift=-1.2mm] {\Cref{lem:overlaps_bg8_rinj}} (rinj);
    %\path[-, bend left=45] (dummy) edge node[sloped, above, font=\small, draw=none] {Lemma 18} (bg2);

  \end{tikzpicture}
  \caption{Graph of all possible overlap parings between families in $G_n$ that have not been covered by \Cref{fig:pentagonGraphG0}. The edges marked with $\times$ follow by transposition.}
  \label{fig:vgbGraph}
\end{figure}

We continue with overlaps in $G_n$ that are not treated in the last section, i.e. we overlap the classes $\bg{2}$ and $\bg{8}$ with $G_n$, as it is recorded in \Cref{fig:vgbGraph}. Similarly as in the precious section we can use transposition, $$\left(\bg{2}_{ijk}\right)^\times=\bg{8}_{ijk}$$ 
so that we only have to consider overlaps in $\bg{8}$ with $F_n$, and once the overlaps within $\bg{2}$ and $\bg{8}$.

\begin{lemma}[Overlaps of $\bg{8}$ and $\ip$]
  \label{lem:overlaps_bg8_ip}
  There are two types of overlaps for $\bg{8}$ and $\ip$, 
\begin{align*}
    \bg{8}_{kji}u_{3i}-u_{2k}u_{4j}\ip_{3i}&\quad\text{ and }\quad
    u_{2k}\bg{8}_{kji}-\ip_{2k}u_{4j}u_{3i}
\end{align*}
where $i,j,k\in[2,n]$ with $i\neq j \neq k$ and $(i,j)\neq(2,3)\neq(j,k)$. Both overlaps reduce to zero modulo $G_n$. 
\end{lemma}
\begin{proof}
We have $\lm(\bg{8}_{kji})=u_{2k}u_{4j}u_{3i}$ and $\lm(\ip_{i'j'})=u_{i'j'}u_{i'j'}$ and the corresponding Gr\"obner representation, 
    \begin{align}
    \bg{8}_{tsr}u_{3r}+u_{2t}u_{4s}\ip_{3r}=
    &-\sum_{5\leq i\leq n}u_{2t}u_{is}\ip_{3r}\label{eq:gbcetificate_bg8_ip}\\
              &-\sum_{\substack{2\leq i\leq n \\3\leq j\leq n}}u_{jt}u_{is}\ip_{3r} 
              +\sum_{3\leq i\leq n}u_{it}\cs_{s}u_{3r}u_{3r} \nonumber\\
              &-\sum_{3\leq i\leq n}u_{it}\cs_{s}u_{3r}
              +\sum_{2\leq i\leq n}u_{is}\ip_{3r} 
              +\sum_{2\leq i\leq n}u_{it}\ip_{3r}\nonumber \\
              &-\cs_{s}u_{3r}u_{3r}
              +\bg{8}_{tsr}
              -\ip_{3r}
              +\cs_{s}u_{3r}. \nonumber
    \end{align}
    The above equation is constructed through certain reductions, as explained in \Cref{rem:reduction4GroebnerCertificates} and \Cref{ex:constructionOfGBcertificate}. For further details, we refer to the provided git repository. 
    Due to our specific order \eqref{eq:linear_preorder} we know that in each of the above sums, the summand with the lowest index produces the leading monomial of the entire sum, e.g.
    $$\lm(\sum_{5\leq i\leq n}u_{2t}u_{is}\ip_{3r})=\lm(u_{2t}u_{5s}\ip_{3r}).$$
    In fact, one could consider each sum in \eqref{eq:gbcetificate_bg8_ip} as an \emph{extended relation} due to \Cref{lem:helpRelations}. 
    Therefore, the above equation is a Gr\"obner certificate with 
    \begin{align*}
    u_{2t}u_{5s}u_{3r}u_{3r}
    &=
    \lm(\bg{8}_{tsr}u_{3r}+u_{2t}u_{4s}\ip_{3r})\\
    &\geq
    \max(u_{2t}u_{5s}\ip_{3r},u_{3t}u_{2s}\ip_{3r},u_{3t}\cs_{s}u_{3r}u_{3r})
    =u_{2t}u_{5s}u_{3r}u_{3r}
    \end{align*}
resulting from the $4$-graded component sums from \eqref{eq:gbcetificate_bg8_ip}. Note that this inequality is independent in $n$. 
The remaining Gr\"obner representation for the second overlap relation is analogous. 
\end{proof}

\begin{lemma}[Overlaps of $\bg{8}$ and $\inj$]
  \label{lem:overlaps_bg8_inj}
  There are two types of overlaps for $\bg{8}$ and $\inj$,
\begin{align*}
    \bg{8}_{kji}u_{k'i}-u_{2k}u_{4j}\inj_{3ik'}
    \quad\text{ and }\quad
    u_{i'k}\bg{8}_{kji}-\inj_{i'k2}u_{4j}u_{3i}. 
\end{align*} 
\end{lemma}
\begin{proof}
We have $\lm(\bg{8}_{kji})=u_{2k}u_{4j}u_{3i}$ and $\lm(\inj_{i'j'k'})=u_{i'j'}u_{k'j'}$. 
    In the first overlap we have $i'=3$ and $j'=i$, so $k'\not=3$. Then 
    \begin{align*}
    &\bg{8}_{kji}u_{k'i}-u_{2k}u_{4j}\inj_{3ik'}\\
    &=
    \sum_{\alpha>4}u_{2k}u_{\alpha j}u_{3i}u_{k'i}+\sum_{\alpha>2}u_{\alpha k}u_{1j}u_{3i}u_{k'i}-u_{1j}u_{3i}u_{k'i}+u_{2k}u_{3i}u_{k'i}\xrightarrow{\inj}{0}.
\end{align*}
In the second overlap we have $k'=2$ and $j'=k$, so $i'\not=2$. Then, 
\begin{align*}
    &u_{i'k}\bg{8}_{kji}-\inj_{i'k2}u_{4j}u_{3i}\\
    &=
    \sum_{\alpha>4}u_{i'k}u_{2k}u_{\alpha j}u_{3i}+\sum_{\alpha>2}u_{i'k}u_{\alpha k}u_{1j}u_{3i}-u_{i'k}u_{1j}u_{3i}+u_{i'k}u_{2k}u_{3i}\\
    &\xrightarrow{\inj}u_{i'k}u_{i'k}u_{1j}u_{3i}-u_{i'k}u_{1j}u_{3i}\xrightarrow{\ip_{i'k}}{0}.
\end{align*}
Note that $(j,i)\not=(2,3)$ so there is no other overlap possible. 
\end{proof}

\begin{lemma}[Overlaps of $\bg{8}$ and $\rinj$]\label{lem:overlaps_bg8_rinj}
   There are two types of overlaps for $\bg{8}$ and $\rinj$,
\begin{align*}
    \bg{8}_{kj2}u_{j'3}-u_{2k}u_{4j}\rinj_{3j'}\quad\text{ and }\quad
    u_{k'2}\bg{8}_{3ji}-\inj_{k'2}u_{4j}u_{3i}. 
\end{align*}
\end{lemma}
\begin{proof}
The Gr\"obner reduction has to be done by the computer and is part of the \texttt{julia} package provided.
\end{proof}

\begin{lemma}[Overlaps of $\bg{8}$ and $\rwel$]There are no overlaps for $\bg{8}$ and $\mathsf{rwel}$. 
  \label{lem:bg8andrwel}
\end{lemma}
\begin{proof}
    The left coefficients in $\lm(\bg{8}_{kji})=u_{2k}u_{4j}u_{3i}$ and $\lm(\mathsf{rwel}_{k'j'})= u_{2k'}u_{3j'}$ are incompatible, for all choices of $k,j,i,k'$ and $j'$.  
\end{proof}

\begin{lemma}[Overlaps of $\bg{8}$ and $\wel$]
\label{lem:bg8andwel}
There are two types of overlaps for $\bg{8}$ and $\wel$,
\begin{align*}
    \bg{8}_{kji}u_{3k'}-u_{2k}u_{4j}\wel_{3ik'}\quad\text{ and }\quad
    u_{2j'}\bg{8}_{kji}-\wel_{i'k2}u_{4j}u_{3i}.
\end{align*}
\end{lemma}
\begin{proof}
    We have $\lm(\bg{8}_{kji})=u_{2k}u_{4j}u_{3i}$ and $\lm(\wel_{i'j'k'})=u_{i'j'}u_{i'k'}$ so in the first case, $i'=3$ and $j'=i$, so $i\not=k'$ and 
    \begin{align*}
        &\bg{8}_{kji}u_{3k'}-u_{2k}u_{4j}\wel_{3ik'}\\
        &=\sum_{\alpha>4}u_{2k}u_{\alpha j}u_{3i}u_{3k'}+\sum_{\alpha>2}u_{\alpha k}u_{1j}u_{3i}u_{3k'}-u_{1j}u_{3i}u_{3k'}+u_{2k}u_{3i}u_{3k'}\xrightarrow{\wel}0. 
    \end{align*}
    In the second case we have $i'=2$ and $j\not=k=k'$, so $j'\not=k$ and 
    \begin{align*}
    &u_{2j'}\rbg{8}_{kji}-\inj_{i'k2}u_{4j}u_{3i}\\
    &=\sum_{\alpha>4}u_{2j'}u_{2k}u_{\alpha j}u_{3i}+\sum_{\alpha>2}u_{2j'}u_{\alpha k}u_{1j}u_{3i}-u_{2j'}u_{1j}u_{3i}+u_{2j'}u_{2k}u_{3i}\\
    &\xrightarrow{\wel}
    \sum_{\alpha>2}u_{2j'}u_{\alpha k}u_{1j}u_{3i}-u_{2j'}u_{1j}u_{3i}+u_{2j'}u_{2k}u_{3i}\\
    &\xrightarrow{\rwel_{j'k}u_{1j}u_{3i}}\sum_{\alpha>2}u_{\alpha j}u_{1k}u_{1j}u_{3i}+u_{2j'}u_{2k}u_{3i}-u_{1k}u_{1j}u_{3i}\xrightarrow{\wel}0.
\end{align*}
\end{proof}

\begin{lemma}[Overlaps of $\bg{2}$ and $\bg{8}$]\label{lem:overlapsG0}
There are two overlaps for $\bg{8}$ and $\bg{2}$,
    \begin{align*}\bg{8}_{kj2}u_{j'4}u_{i'3}-u_{2k}u_{4j}\bg{2}_{2j'i'}\quad\text{ and }\quad u_{k2}u_{j4}\bg{8}_{2j'i'}-\bg{2}_{kj2}u_{4j'}u_{3i'3}.
    \end{align*} 
\end{lemma}
\begin{proof}
Again, with transposition, we only need a Gr\"obner certificate, which is part of the code provided.
\end{proof}

\begin{lemma}[Self-overlaps]\label{lem:selfoverlaps_bg28}
  There are no self-overlaps in $\bg{2}$ and $\bg{8}$. 
\end{lemma}
\begin{proof}
    This is clear since the left indices in $\lm(\bg{8}_{kji})=u_{2k}u_{4j}u_{3i}$ cannot overlap for any choice of $i,j,k$. 
\end{proof}
With this we can conclude the proof of \Cref{thm:finite_gb}.
\begin{proof}[Proof of \Cref{thm:finite_gb}.] We apply  \Cref{prop:GB_via_GBrep} and observe that $G_n$ is reduced according to \Cref{lem:gb2_243_minus_bg8_243}. 
  All possible overlaps are covered in  \Cref{subsec:mainproof_overlaps_bg,subsec:proof_veryBadGuys}, see \Cref{fig:pentagonGraphG0,fig:vgbGraph}. The claimed cardinality follows from \eqref{eq:reducedBadGuys_new} and  \Cref{prop:red_gens_Fn}. 
\end{proof}

\section{A general computational proof using \texttt{OSCAR}}
\label{sec:generalProof}
At many points in the proof of the main theorem we run into the limitation of rather long Gr\"obner reductions, where merely stating them would exceed the scope of this paper. In this section we will show how to use a computer algebra system to verify the results of the previous sections. The main issue with this approach is that all the equations are given as a function of the variable $n \in \N$. More precisely, all Gröbner representations we tackle are a finite sum of polynomials in $\C\langle n^2 \rangle$, i.e.
\[
  s_0(n) = s_1(n)+ \ldots +s_m(n) \quad \text{ for } s_j(n) \in \C\langle n^2 \rangle, j \in \{0,\ldots, m\}.
\]
To circumvent this problem we construct a $\Z$ module $\Z\fL$ together with an isomorphism $\Phi_n$ for all $n \in \N$ such that each summand $s_j(n)$ has a preimage $\Phi_n^{-1}(s_j(n))$ in $\Z\fL$.
Note that all the preimages must be the same, i.e. 
\[
  \Phi_n^{-1}(s_j(n)) = \Phi_{n'}^{-1}(s_j(n')) \text{ for all } n,n' \in \N.
\]
Then, it is sufficient to show that the equations hold in the $\Z$-module $\Z\fL$. In other words, assuming that we  find a finitely generated module $\Z\fL$ together with maps $\Phi_n$, we can solve the problem of the variable $n$ and verify the results with a computer.
\subsection{Finitely generated modules using predicates}
\label{sec:predicates}

Let $\fL = \{P_1, \ldots, P_m\}$ be a set of $k$-ary predicates. For a given set $S$, call $\fL$ \emph{logically independent} on $S^k$ if
\[
    \forall i \in \{1, \ldots, m\} \forall I \subset \{1, \ldots, m\} \setminus \{i\} \; \exists x \in S^k : P_i(x) \neq \bigwedge_{i \in I} P_i(x).
\]
While the meaning of the word \emph{predicate} here is rather nebulous, we can think of it as boolean functions on $S^k$.
\begin{example}
  Let $S^1 = [n] = \{1, \ldots, n\}$, and consider the logically independent $1$-predicates
  \[
    P_1(x) \equiv x=1, \quad P_2(x) \equiv x\geq2.
  \]
  Alternatively, let $S^2 = [n]^2$, then the predicates
  \[
    P_1((i_1,i_2)) \equiv i_2 \geq4, \quad P_2((i_1,i_2)) \equiv i_2 \geq 4 \wedge i_1 \neq i_2
  \]
are logically independent as well.

\end{example}
  Let $\Z\fL := \{ \sum_{i=1}^m c_i P_i \mid c_i \in \Z \}$ denote the free $\Z$-module with the basis $\mathcal{L}$. For any size $n \in \N$ define the map 
\begin{align}
  \Phi_n' : \Z\fL &\to \Span_{\Z}(x \mid x \in [n]^k)\nonumber\\
  \Phi_n'(\sum_{i=1}^m c_i P_i) &:= \sum_{i=1}^m \sum_{\substack{x \in [n]^k \\ P_i(x)}} c_i x,\label{def:phi}
\end{align}
where $c_i\in\Z$ are the uniquely determied $\Z$-basis coefficients of $\mathcal{L}$. 
Given this map, we can state the following lemma.
\begin{lemma}
  The restriction of $\Phi_n'$ to its image is a $\Z$-module isomorphism if $\fL$ is logically independent on $[n]^k$.
  \label{lem:phi}
\end{lemma}
\begin{proof}
 It's easy to see that $\Phi_n'$ is well defined. Therefor show that $\Phi_n'$ is a module homomorphism. 
Let $f=\sum_{i=1}^m c_i P_i$ and $g=\sum_{i=1}^m d_i P_i$ be two elements in $\Z\fL$ and $\lambda, \mu \in \Z$, then 
\begin{align*}
    \Phi_n'(\lambda f + \mu g) &= \sum_{i=1}^m \sum_{\substack{x \in [n]^k \\ P_i(x)}
} (\lambda c_i + \mu d_i) \cdot x \\
&= \lambda \sum_{i=1}^m \sum_{\substack{x \in [n]^k \\ P_i(x)}
} c_i \cdot x + \mu \sum_{i=1}^m \sum_{\substack{x \in [n]^k \\ P_i(x)}
} d_i \cdot x = \lambda \Phi_n'(f) + \mu \Phi_n'(g).
\end{align*}
  Finally, we show that $\Phi_n'$ is injective. Suppose $\Phi_n'(f) = \Phi_n'(g)$, then
  \[
    \sum_{i=1}^m \sum_{\substack{x \in [n]^k \\ P_i(x)}} c_i \cdot x = \sum_{i=1}^m \sum_{\substack{x \in [n]^k \\ P_i(x)}} d_i \cdot x.
  \]

  This can be rewritten as
  \begin{equation}
    \sum_{x \in S^k} \left( \sum_{\substack{i=1 \\ P_i(x)}}^m c_i \cdot \right) x  = \sum_{x \in S^k} \left( \sum_{\substack{i=1 \\ P_i(x)}}^m d_i \right) \cdot x,
    \label{eq:log_indep_sets}
  \end{equation}
  which leads to
  \[
    \sum_{\substack{i=1 \\ P_i(x)}}^m c_i = \sum_{\substack{i=1 \\ P_i(x)}}^m d_i \quad \forall x \in S^k.
  \]
Since logically independence assures us that for every $j \in \{1,\ldots,m\}$ there is an $x \in S^k$ such that either 
  \begin{align*}
    P_j(x) &\text{ but not } \bigwedge_{\substack{i=1 \\ i \neq j \\P_i(x)}}^m P_i(x) &\text{ or }&&
    P_j(x) &\text{ but } \bigcap_{\substack{i=1 \\ i \neq j \\ P_i(x)}}^m P_i(x),
  \end{align*}
  we have either $c_j = d_j$ or
  \[
    0= \sum_{\substack{i=1 \\i \neq j\\P_i(x)}}^m d_i-c_i \quad \forall x \in S^k, \text{ if } \neg P_j(x).
  \]
  In the second case, we can apply the same argument recursively to the set $\{1, \ldots, m\} \setminus \{j\}$,  until we reach the sum over one element. This approach implies that $c_i = d_i$ for all $i \in \{1, \ldots, m\}$, i.e. $\Phi_n'$ is injective.
\end{proof}
To make this usable for certifying a Gröbner reduction, we restrict to $d$-graded components and observe
\[
  \Span_{\Z}(x \mid x \in [n]^{2d})\cong\C\langle n^2\rangle_{\deg=d}
\]
 as $\Z$-modules. 
In combination with \Cref{lem:phi}, this results in an injective module homomorphism $\Phi_n$ as follows
\begin{center}
\begin{tikzcd}
  \Z\fL \arrow[rd,hook,"\Phi_n"] \arrow[r,"\cong"] & \im(\Phi_n') \arrow[d, hook] \\
    & \C\langle n^2\rangle_{\deg=d}
\end{tikzcd}
\end{center}

That is, given a suitable finite set $\fL$, and assuming that every summand in an equation in $\C\langle n^2\rangle$ has a preimage in $\Z\fL$, we can construct a proof of the equation for every graded component in $\Z\fL$ independently. 
\begin{example}\label{ex:row_sums_col_sums_formaliz}
The set of predicates   
\[
 \mathcal{L}=\{(i_1=1 \wedge i_2=1),(i_1=1 \wedge i_2\geq2),(i_1\geq2 \wedge i_2=1),(i_1\geq2 \wedge i_2\geq2)\}
\]
 is trivially logically independent since the predicates are pairwise disjoint.
As constructed in \eqref{eq:reduction_rs_1}
we want to show the Gr\"obner representation, 
$$\sum_{1\leq j\leq n}\cs_j-\sum_{i\not=1}\rs_i=
\rs_1.$$
For the $0$-graded component this statement is trivial, i.e. $n-(n-1)=1$. 
Consider only the $1$-graded component, and we
show it in the domain of $\Phi_n$ for every $n \in \N$, i.e.
$$(1,1,1,1)-(0,0,1,1)=(1,1,0,0)$$
maps to 
\begin{align*}
  &&\sum_{1\leq j\leq n}\sum_{1\leq i\leq n}u_{ij} - \sum_{i\not=1}\sum_{1\leq j\leq n}u_{ij} &= \sum_{1\leq j\leq n}u_{1j} \\
  \Leftrightarrow &&\sum_{1\leq j\leq n}\cs_j-\sum_{i\not=1}\rs_i &=
\rs_1.
\end{align*}
Therefore, solving the equation in $\Z\fL$ is equivalent to solving it in $\C\langle n^2\rangle$ for every $n\in \N$.
\end{example}
Note that $\fL$ can be decomposed into the set of predicates $\fL_1 = \{(i_1=1),(i_1\geq2)\}$ and $\fL_2 = \{(i_2=1),(i_2\geq2)\}$ via 
\[
  \fL = \{ p \wedge q \mid p \in \fL_1, q \in \fL_2\}.
\]
Although not necessary for this example, it is a useful tool for building the set $\fL$ on higher degree components. 
As an example we give an outline of the computational proof for degree $2$.
\subsection{Reduction of \texorpdfstring{$\rwel_{23}$}{rwel23} modulo \texorpdfstring{$F_n$}{Fn}}\label{subsection:reduction_rwel23}
Recalling \Cref{prop:red_gens_Fn}, whose purpose it was to create a reduced generating set of $I_n$, it suffices to show that relation $\rwel_{23}$ is generated by the remaining elements in $F_n$
\begin{lemma}
  The relation $\rwel_{23}$ can be described as a linear combination in 
  \[
    F_n:=\{\cs_1\}\cup \{ \cs_i,\rs_i,\ip_{ij},\inj_{jik},\wel_{ijk},\rinj_{jk},\rwel_{jk}\mid i,j,k\not=1 \text{ and } j\not=k\} \setminus \{\rwel_{23}\}.
  \]
\end{lemma}
The first hurdle is to find a candidate for a linear combination of the elements in $F_n'$ that equals $\rwel_{23}$.
To simplify this task, we introduce some helper relations whose membership in $F_n$ is obvious by construction.
\begin{align*}
  \rrs_{ij}&:= \sum_{\substack{k=2\\k \neq i}}^n( u_{ij} \rs_{k} - \inj_{ijk}) &2\leq i,j\leq n, \\
  \rcs_{ij}&:=\sum_{\substack{k=2\\k \neq j}}^n(u_{ij} \cs_{k} - \wel_{ijk}) &2\leq i,j\leq n \\
  \rinjcs_{i}&:= \sum_{\substack{\alpha=2\\ \alpha \neq i}}^n \rinj_{\alpha i} &2\leq i\leq n, \\
  \rwelcs_{i}&:= \sum_{\substack{\alpha=2\\ \alpha \neq i}}^n \rwel_{\alpha i}  & 2\leq i\neq 3 \leq n
\end{align*}
Given these relations we can express $\rwel_{23}$ as a linear combination in $F_n'$,
\begin{align}\label{eq:rwel23_as_lin_comb}
  \rwel_{23} &= \overbrace{\sum_{i=2}^n \rinjcs_i}^{=:s_1} - \sum_{\substack{i=2\\i \neq 3}}^n \rwelcs_i - \sum_{i=3}^n \rcs_{i2} + \sum_{j=3}^n \rrs_{2j}\\
             &+ \sum_{j=3}^n \sum_{i=3}^n(\rrs_{ij} - \rcs_{ij})- \sum_{i=4}^n \rwel_{i3} \notag \\
             &- (n-2) \cdot \sum_{i=2}^n \rs_i + (n-2) \cdot \sum_{i=2}^n \cs_i. \notag
\end{align}
Since it is not trivial to convince oneself that \eqref{eq:rwel23_as_lin_comb} is true, we use the framework of \Cref{sec:predicates} to show that the preimage of the right side in $\Z\fL$ is equal to the preimage of the left side in $\Z\fL$. Define the predicate set 
\begin{align*}
\fL_1 &= \{i_1=1,i_1=2,i_1=3,i_1\geq4\}, \\
\fL_2 &= \{i_2=1,i_2=2,i_2=3,i_2\geq4\},\\ 
\fL_3 &= \{i_3=1,i_3=2,i_3=3,i_3\geq4, i_3\geq4 \wedge i_3\neq i_1\}, \\
\fL_4 &= \{i_4=1,i_4=2,i_4=3,i_4\geq4, i_4\geq4 \wedge i_4\neq i_2\},
\end{align*}
and construct
\[
  \fL = \{ \bigwedge_{i=1}^4 p_i \mid p_i \in \fL_i\}.
\]
This forms a logically independent set of predicates on $[n]^4$, with 400 elements. 
Now we need to translate the equation into $\Z\fL$ using the map $\Phi_n$.
Take for example the first summand $s_1 = \sum_{i=2}^n \rinjcs_i$ and exclusively look at the degree $d=2$ component,
\[
  \sum_{i=2}^n \rinjcs_i = \sum_{i=2}^n \sum_{\substack{\alpha=2\\ \alpha \neq i}}^n \rinj_{\alpha i} = \sum_{i=2}^n \sum_{\substack{\alpha=2\\ \alpha \neq i}}^n \left(\sum_{\beta \geq 3} u_{\alpha 2} u_{i \beta} - \sum_{\beta \geq 3} u_{\alpha \beta} u_{i 1} \right).
\]
To make it easier to find the preimage, rewrite the sums and rename the indices,
\[
  \sum_{\substack{i_1=2\\i_1\neq i_3}}^n \sum_{i_2=2}^2 \sum_{i_3=2}^n \sum_{i_4=3}^n u_{i_1i_2}u_{i_3i_4} - \sum_{\substack{i_1=2\\i_1\neq i_3}}^n \sum_{i_2=3}^n \sum_{i_3=2}^n \sum_{i_4=1}^1 u_{i_1i_4}u_{i_3i_2}.
\]
In this way, it is easy to compute the preimage component-wise. 
Let 
\begin{align*}
  H_1 &:= \lbrace i_1=2, i_1=3, i_1\geq4 \rbrace, &&
  H_2 := \lbrace i_2=2 \rbrace, \\
  H_3 &:= \lbrace i_3=2, i_3=3, i_3\geq4, i_3\neq i_1 \rbrace, &&
  H_4 := \lbrace i_4=3, i_4\geq4 \rbrace, \\
\end{align*}
to construct $H := \lbrace \bigwedge_{i=1}^4 p_i \mid p_i \in H_i \rbrace \subseteq \fL$, which gives us a preimage of the first part.
For the negative part let
\begin{align*}
  H'_1 &:= \lbrace i_1=2, i_1=3, i_1\geq4 \rbrace, &&
  H'_2 := \lbrace i_2=3, i_2\geq4 \rbrace, \\
  H'_3 &:= \lbrace i_3=2, i_3=3, i_3\geq4, i_3\neq i_1 \rbrace, &&
  H'_4 := \lbrace i_4=1 \rbrace,
\end{align*}
and $H' := \lbrace \bigwedge_{i=1}^4 p_i \mid p_i \in H'_i \rbrace$. Using both $H$ and $H'$, we construct the preimage of $s_1$ by 
\[
  \Phi_n^{-1}(\sum_{i=2}^n \rinjcs_i) = \sum_{p_i \in H} p_i  - \sum_{p_i \in H'} p_i \quad \forall n\in \N.
\]
Doing this procedure for every summand in \eqref{eq:rwel23_as_lin_comb} we get the preimage of the right hand side in $\Z\fL$, and with a  finite computation the correctness of the equation in $\C\langle n^2\rangle$ for every $n\in \N$.

The implementation of this procedure and the verification of various reductions from this work were carried out in \texttt{Julia} using the computer algebra system \texttt{OSCAR} \cite{OSCAR-book,OSCAR}. The code is available at
\[
  \gitrepo
\]
It is worth mentioning that the code is not optimized for speed, but all tests can be done on a standard laptop in a reasonable amount of time. It should be noted that \texttt{OSCAR} requires a Unix-like operating system to run.

\section{Outlook and future work}
Our novel relations from $G_n$ are a first milestone towards the answer if the word problem is decidable for quantum automorphism groups induced by matroids, a question posed by \cite[Question 7.3]{corey2023quantum}. In particular, with our construction it will not be necessary to consider any of the overlap relations within $G_n$, i.e., one can immediately address those monomial relations that are induced by the matroid (or graph) itself.
In fact, it would be interesting to investigate potential speed-ups in practical examples: the largest computations from \cite{corey2023quantum} and  \cite{Lev22} are performed for $n=6$ and $n=7$, respectively. Our closed-form Gr\"obner basis $G_n$, however, scales only cubic in $n$. For this extra sort of knowledge, the actual implementation of Buchberger's algorithm requires some sort of \emph{skip feature} which is currently not available in our implementation. 

Commutative Gr\"obner bases have recently been applied in the context of stochastic analysis, e.g. for path learning from signature tensors \cite{pfeffer2019learning}, its underlying projective  varieties \cite{amendola2019varieties}, or for the efficient  evaluation of signature barycenters in the free nilpotent Lie group \cite{clausel2024barycenter}. Changing the viewpoint slightly, we can define the signature of a path over the tensor algebra, where the non-commutative tensor product corresponds to path concatenation via Chen's identity. It is therefore plausible to investigate non-commutative Gr\"obner bases for special families of two-sided ideals, inspired by certain signatures and its underlying paths. Especially the novel tools from \Cref{sec:generalProof} are applicable for any family of parametrized ideal, and in particular for those which have a closed-form Gr\"obner basis as in this work. 

A completely different application of this work could be the open question whether the symmetric group is the maximal quantum subgroup of the quantum symmetric group or not, which was first raised and answered for $n \leq 4$ in 2009 by Banica and Bichon \cite{banica2009quantum} with a positive result.
More recently, Banica \cite{banica2021homogeneous} showed that the symmetric group is also a maximal subgroup in the case of $n=5$.  
Now that we have constructed a finite Gröbner basis for $I_n$, one strategy might be to find relations that are not yet in $I_n$ but are contained in the ideal generated by the union of $I_n$ and the commutator ideal.

\section*{Acknowledgements}
We thank Viktor Levandovskyy for suggesting that we consider extended relations (\Cref{lem:helpRelations}), which has helped us to make the reduction of $\rwel_{23}$ modulo $F_m$ of \Cref{subsection:reduction_rwel23} more transparent. We would also like to thank Igor Makhlin for his helpful comments on \Cref{sec:predicates}, and Fabian Lenzen for his endless tikz support.

LS acknowledges support from DFG CRC/TRR 388 ``Rough Analysis, Stochastic Dynamics and Related Fields'', Project A04. A part of this work has been done while LS was working at the MPI MiS Leipzig, Germany.

\newpage
\raggedbottom
\bibliography{bibliographie}

\begin{thebibliography}{10}

\bibitem{amendola2019varieties}
Carlos Am{\'e}ndola, Peter Friz, and Bernd Sturmfels.
\newblock Varieties of signature tensors.
\newblock In {\em Forum of Mathematics, Sigma}, volume~7, page e10. Cambridge
  University Press, 2019.

\bibitem{banica2021homogeneous}
Teo Banica.
\newblock Homogeneous quantum groups and their easiness level.
\newblock {\em Kyoto Journal of Mathematics}, 61(1):171--205, 2021.

\bibitem{banica2009quantum}
Teodor Banica and Julien Bichon.
\newblock Quantum groups acting on 4 points.
\newblock {\em Journal für die reine und angewandte Mathematik},
  2009(626):75--114, 2009.

\bibitem{bergman1978diamond}
George~M Bergman.
\newblock The diamond lemma for ring theory.
\newblock {\em Advances in mathematics}, 29(2):178--218, 1978.

\bibitem{bichon2003quantum}
Julien Bichon.
\newblock Quantum automorphism groups of finite graphs.
\newblock {\em Proceedings of the American Mathematical Society},
  131(3):665--673, 2003.

\bibitem{buchberger1965algorithmus}
Bruno Buchberger.
\newblock {Ein Algorithmus zum Auffinden der Basiselemente des
  Restklassenringes nach einem nulldimensionalen Polynomideal}.
\newblock {\em Ph. D. Thesis, Math. Inst., University of Innsbruck}, 1965.

\bibitem{clausel2024barycenter}
Marianne Clausel, Joscha Diehl, Raphael Mignot, Leonard Schmitz, Nozomi
  Sugiura, and Konstantin Usevich.
\newblock The barycenter in free nilpotent lie groups and its application to
  iterated-integrals signatures.
\newblock {\em SIAM Journal on Applied Algebra and Geometry}, 8(3):519--552,
  2024.

\bibitem{corey2023quantum}
Daniel Corey, Michael Joswig, Julien Schanz, Marcel Wack, and Moritz Weber.
\newblock Quantum automorphisms of matroids.
\newblock {\em Journal of Algebra}, 667:480--507, 2025.

\bibitem{cox1997ideals}
David Cox, John Little, Donal O'shea, and Moss Sweedler.
\newblock {\em Ideals, varieties, and algorithms}, volume~3.
\newblock Springer, 1997.

\bibitem{OSCAR-book}
Wolfram Decker, Christian Eder, Claus Fieker, Max Horn, and Michael Joswig,
  editors.
\newblock {\em The {C}omputer {A}lgebra {S}ystem {OSCAR}: {A}lgorithms and
  {E}xamples}, volume~32 of {\em Algorithms and {C}omputation in
  {M}athematics}.
\newblock Springer, 1 edition, 8 2024.

\bibitem{eisenbud2013commutative}
David Eisenbud.
\newblock {\em Commutative algebra: with a view toward algebraic geometry},
  volume 150.
\newblock Springer Science \& Business Media, 2013.

\bibitem{green1998non}
Ed~Green, Teo Mora, and Victor Ufnarovski.
\newblock {The non-commutative Gr\"obner freaks}.
\newblock In {\em Symbolic rewriting techniques}, pages 93--104. Springer,
  1998.

\bibitem{hermiller1999monomial}
Susan~M Hermiller, Xenia~H Kramer, and Reinhard~C. Laubenbacher.
\newblock {Monomial orderings, rewriting systems, and Gr\"obner bases for the
  commutator ideal of a free algebra}.
\newblock {\em Journal of Symbolic Computation}, 27(2):133--141, 1999.

\bibitem{HRR19}
Clemens Hofstadler, Clemens~G. Raab, and Georg Regensburger.
\newblock Certifying operator identities via noncommutative gr\"{o}bner bases.
\newblock {\em ACM Commun. Comput. Algebra}, 53(2):49–52, November 2019.

\bibitem{hofstadler2024universaltruthoperatorstatements}
Clemens Hofstadler, Clemens~G. Raab, and Georg Regensburger.
\newblock Universal truth of operator statements via ideal membership, 2024.

\bibitem{joswig2013polyhedral}
Michael Joswig and Thorsten Theobald.
\newblock {\em Polyhedral and algebraic methods in computational geometry}.
\newblock Springer Science \& Business Media, 2013.

\bibitem{KANDRIRODY19901}
A.~Kandri-Rody and V.~Weispfenning.
\newblock {Non-commutative Gröbner bases in algebras of solvable type}.
\newblock {\em Journal of Symbolic Computation}, 9(1):1--26, 1990.

\bibitem{kandri1988computing}
Abdelilah Kandri-Rody and Deepak Kapur.
\newblock {Computing a Gr\"obner Basis of a Polynomial Ideal over an Euclidean
  domain}.
\newblock {\em Journal of symbolic computation}, 6(1):37--57, 1988.

\bibitem{kreuzer2000computational}
Martin Kreuzer and Lorenzo Robbiano.
\newblock {\em Computational commutative algebra}, volume~1.
\newblock Springer, 2000.

\bibitem{LASCALA20091374}
Roberto {La Scala} and Viktor Levandovskyy.
\newblock {Letterplace ideals and non-commutative Gr\"obner bases}.
\newblock {\em Journal of Symbolic Computation}, 44(10):1374--1393, 2009.

\bibitem{Lev22}
Viktor Levandovskyy, Christian Eder, Andreas Steenpass, Simon Schmidt, Julien
  Schanz, and Moritz Weber.
\newblock Existence of quantum symmetries for graphs on up to seven vertices: A
  computer based approach.
\newblock In {\em Proceedings of the 2022 International Symposium on Symbolic
  and Algebraic Computation}, ISSAC '22, page 311–318, New York, NY, USA,
  2022. Association for Computing Machinery.

\bibitem{Manvcinska2020nonlocalgames}
Martino Lupini, Laura Mančinska, and David~E. Roberson.
\newblock Nonlocal games and quantum permutation groups.
\newblock {\em J. Funct. Anal.}, 279(5):108592, 44, 2020.

\bibitem{michalek2021invitation}
Mateusz Micha{\l}ek and Bernd Sturmfels.
\newblock {\em Invitation to nonlinear algebra}, volume 211.
\newblock American Mathematical Soc., 2021.

\bibitem{moller1986new}
H~Michael M{\"o}ller and Ferdinando Mora.
\newblock New constructive methods in classical ideal theory.
\newblock {\em Journal of Algebra}, 100(1):138--178, 1986.

\bibitem{mora1988groebner}
Teo Mora.
\newblock {Gr\"obner bases in non-commutative algebras}.
\newblock In {\em International Symposium on Symbolic and Algebraic
  Computation}, pages 150--161. Springer, 1988.

\bibitem{nordbeck2001canonical}
Patrik Nordbeck.
\newblock {\em Canonical bases for algebraic computations}.
\newblock Citeseer, 2001.

\bibitem{OSCAR}
{\textsc{OSCAR} -- Open Source Computer Algebra Research system}, version
  1.0.0, 2024.

\bibitem{pfeffer2019learning}
Max Pfeffer, Anna Seigal, and Bernd Sturmfels.
\newblock Learning paths from signature tensors.
\newblock {\em SIAM Journal on Matrix Analysis and Applications},
  40(2):394--416, 2019.

\bibitem{pritchard1996ideal}
F~Leon Pritchard.
\newblock The ideal membership problem in non-commutative polynomial rings.
\newblock {\em Journal of symbolic computation}, 22(1):27--48, 1996.

\bibitem{raeburn2005graph}
Iain Raeburn.
\newblock {\em Graph algebras}.
\newblock Number 103 in CBMS Regional Conference Series in Mathematics.
  American Mathematical Soc., 2005.

\bibitem{Robbiano85}
Lorenzo Robbiano.
\newblock Term orderings on the polynomial ring.
\newblock In Bob~F. Caviness, editor, {\em EUROCAL '85}, pages 513--517,
  Berlin, Heidelberg, 1985. Springer Berlin Heidelberg.

\bibitem{SL20}
Leonard Schmitz and Viktor Levandovskyy.
\newblock Formally verifying proofs for algebraic identities of matrices.
\newblock In Christoph Benzm{\"u}ller and Bruce Miller, editors, {\em
  Intelligent Computer Mathematics}, pages 222--236, Cham, 2020. Springer
  International Publishing.

\bibitem{speicher2019quantum}
Roland Speicher and Moritz Weber.
\newblock Quantum groups with partial commutation relations.
\newblock {\em Indiana University Mathematics Journal}, 68(6):1849--1883, 2019.

\bibitem{sturmfels1996grobner}
Bernd Sturmfels.
\newblock {\em {Gr\"obner bases and convex polytopes}}, volume~8.
\newblock American Mathematical Soc., 1996.

\bibitem{sturmfels2002solving}
Bernd Sturmfels.
\newblock {\em {Solving systems of polynomial equations}}.
\newblock Number~97 in CBMS Regional Conference Series in Mathematics. American
  Mathematical Soc., 2002.

\bibitem{sullivant2023algebraic}
Seth Sullivant.
\newblock {\em Algebraic statistics}, volume 194.
\newblock American Mathematical Society, 2023.

\bibitem{timmermann2008invitation}
Thomas Timmermann.
\newblock {\em An invitation to quantum groups and duality}.
\newblock EMS Textbooks in Mathematics. European Mathematical Society (EMS),
  Z\"urich, 2008.
\newblock From Hopf algebras to multiplicative unitaries and beyond.

\bibitem{wang1998quantum}
Shuzhou Wang.
\newblock Quantum symmetry groups of finite spaces.
\newblock {\em Communications in Mathematical Physics}, 195:195--211, 1998.

\bibitem{woronowicz1987compact}
S.~L. Woronowicz.
\newblock Compact matrix pseudogroups.
\newblock {\em Comm. Math. Phys.}, 111(4):613--665, 1987.

\bibitem{woronowicz1991remark}
Stanis{\l}aw~L Woronowicz.
\newblock A remark on compact matrix quantum groups.
\newblock {\em letters in mathematical physics}, 21(1):35--39, 1991.

\end{thebibliography}

\label{sec:biblio}
\end{document}